\documentclass[a4paper,12pt]{amsart}

\usepackage{amsmath,amssymb,amsfonts,mathrsfs,amsthm,stmaryrd,mathtools,enumerate}
\usepackage[all]{xy}
\usepackage{url}
\usepackage[colorlinks=true]{hyperref}
\usepackage[nobysame,alphabetic,initials,msc-links]{amsrefs}
\usepackage{cleveref}

\DefineSimpleKey{bib}{arxiv}
\newcommand\myurl[1]{\url{#1}}

\def\MR#1{\quad \href{http://www.ams.org/mathscinet-getitem?mr=#1}{MR#1}}

\DefineSimpleKey{bib}{how}
\renewcommand{\eprint}[1]{#1}
\BibSpec{misc}{%
  +{}{\PrintAuthors}  {author}
  +{,}{ \textit}      {title}
  +{,}{ }             {how}
  +{}{ \parenthesize} {date}
  +{,} { available at \eprint}        {eprint}
  +{,}{ available at \url}{url}
  +{,}{ }             {note}
  +{.}{}              {transition}
}

\newtheorem{thrm}{Theorem}[section]
\newtheorem{prop}[thrm]{Proposition}
\newtheorem{coro}[thrm]{Corollary}
\newtheorem{lemm}[thrm]{Lemma}

\theoremstyle{definition}
\newtheorem{defn}[thrm]{Definition}

\newtheorem{rema}[thrm]{Remark}

\renewcommand{\labelenumi}{(\theenumi)}

\renewcommand{\tilde}[1]{\widetilde{#1}}
\renewcommand{\epsilon}{\varepsilon}
\newcommand{\A}{\mathcal{A}}

\newcommand{\aE}{\acute{E}}
\newcommand{\hE}{\hat{E}}
\newcommand{\tE}{\tilde{E}}
\newcommand{\vE}{\check{E}}
\newcommand{\aF}{\acute{F}}
\newcommand{\hF}{\hat{F}}
\newcommand{\tF}{\tilde{F}}
\newcommand{\vF}{\check{F}}
\newcommand{\aK}{\acute{K}}
\newcommand{\hK}{\hat{K}}
\newcommand{\tK}{\tilde{K}}
\newcommand{\vK}{\check{K}}

\newcommand{\Z}{\mathbb{Z}}
\newcommand{\Q}{\mathbb{Q}}
\newcommand{\R}{\mathbb{R}}
\newcommand{\C}{\mathbb{C}}

\newcommand{\frb}{\mathfrak{b}}
\newcommand{\g}{\mathfrak{g}}
\newcommand{\tg}{{}^{\sim}\g}
\newcommand{\gt}{\g{}^{\sim}}

\newcommand{\hH}{\hat{H}}
\newcommand{\h}{\mathfrak{h}}
\newcommand{\frl}{\mathfrak{l}}

\newcommand{\lt}{\frl_S{}^{\sim}}
\newcommand{\lot}{\frl_0{}^{\sim}}
\newcommand{\tm}{\tilde{m}}
\newcommand{\n}{\mathfrak{n}}

\newcommand{\frp}{\mathfrak{p}}
\newcommand{\frpt}{\frp_S{}^{\sim}}
\newcommand{\pot}{\frp_0{}^{\sim}}
\newcommand{\tfrp}{{}^{\sim}\frp_S}
\newcommand{\fru}{\mathfrak{u}}

\newcommand{\utp}{\fru_S^+{}^{\sim}}
\newcommand{\utm}{\fru_S^-{}^{\sim}}

\newcommand{\eps}{\varepsilon}
\renewcommand{\phi}{\varphi}

\newcommand{\map}[3]{#1\colon#2 \longrightarrow #3}

\newcommand{\abs}[1]{\lvert #1 \rvert}

\newcommand{\cl}[1]{\mathcal{#1}}
\newcommand{\rsfs}[1]{\mathscr{#1}}

\newcommand{\bs}[1]{\boldsymbol{#1}}

\newcommand{\qbin}[3]{\begin{bmatrix} #1 \\ #2 \end{bmatrix}_{#3}}
\newcommand{\tensor}{\otimes}
\newcommand{\htensor}{\hat{\tensor}}
\newcommand{\fin}{\mathrm{f}}

\newcommand{\tend}{\textendash}
\renewcommand{\mod}[1]{#1\text{-}\mathrm{Mod}}
\newcommand{\roots}{R}
\newcommand{\simples}{\Delta}
\newcommand{\bbh}{\llbracket h \rrbracket}

\newcommand{\intsts}{\cl{O}_q^{\A}(B^+B^-)}
\newcommand{\pdqind}[1]{\mathrm{ind}_{\frl_0,q}^{\g,\chi} #1}
\newcommand{\dqind}[1]{\mathrm{ind}_{\frl_S,q}^{\g,\chi} #1}
\newcommand{\qind}[1]{\mathrm{ind}_{\frp_S,q}^{\g} #1}
\newcommand{\deff}{\overset{\text{def}}{\Longleftrightarrow}}

\DeclareMathOperator{\spec}{\mathrm{Spec}}
\DeclareMathOperator{\support}{\mathrm{supp}}

\DeclareMathOperator{\id}{\mathrm{id}}
\DeclareMathOperator{\wt}{\mathrm{wt}}
\DeclareMathOperator{\Hom}{\mathrm{Hom}}
\DeclareMathOperator{\Supp}{\mathrm{supp}}

\DeclareMathOperator{\aspan}{\mathrm{span}}

\DeclareMathOperator{\rep}{\mathrm{Rep}}
\DeclareMathOperator{\modulo}{\mathrm{mod}}

\DeclareMathOperator{\irr}{\mathrm{Irr}}

\title[Poly. families of quantum s.s. coadjoint orbits via deformed QEA]{Polynomial families of quantum semisimple coajoint orbits via deformed quantum enveloping algebras}
\author{Mao Hoshino}
\address{Department of Mathematical Sciences, The University of Tokyo\\
Komaba 3-8-1, Tokyo \mbox{153-8914}, Japan}
\email{mhoshino@ms.u-tokyo.ac.jp}
\thanks{This work was supported by JSPS KAKENHI Grant Number JP23KJ0695 and WINGS-FoPM Program at the University of Tokyo.}

\begin{document}

\begin{abstract}
We construct a polynomial family of semisimple left module categories
over the representation category of the Drinfeld-Jimbo deformation,
with the fusion rule of the representation category of each Levi subalgebra.
In this construction we perform a kind of generalized parabolic induction using a deformed quantum enveloping algebra,
whose definition depends on an arbitrary choice of a positive system and corresponds to
De Commer's definition for the standard positive system. These algebras
define a sheaf of algebras on the toric variety associated to
the root system, which contains the moduli of equivariant Poisson brackets.
This fact finally produces the family of $2$-cocycle.
We also obtain a comparison theorem between our module categories and
module categories induced from our construction for
intermediate Levi subalgebras.
The construction of deformed quantum enveloping algebras and
the comparison theorem are discussed in the integral setting of
Lusztig's sense.
\end{abstract}

\maketitle

\section{Introduction}

Deformation quantization is one of rigorous formulations
of quantization, which naturally arose in quantum physics. 
This physical background is a reason why we typically consider
an algebra of suitable functions on a space with a Poisson structure as
an algebra to be quantized, though it is possible to introduce the notion
of deformation quantization for general Poisson algebras.
Actually finding and classifying deformation quantizations of such function algebras is a natural problem to consider and has been attracting
many physicists and mathematicians.
One of most successful results was accomplished by Kontsevich (\cite{MR2062626}),
who proved the existence of deformation quantizations of each Poisson manifold and their classification theorem.
See \cite{MR1321655} for developments before Kontsevich's work,
including Fedosov's construction
for deformation quantizations of symplectic manifolds (\cite{MR1293654}).

Theory of quantum groups is also a relevant and important
branch of mathematical physics around quantization.
Not surprisingly, deformation quantization deeply relates to it.
For instance a Drinfeld-Jimbo deformation $U_q(\g)$
of a complex semisimple Lie algebra $\g$ is a kind of dual
of a quantization $\cl{O}_q(G)$ of the function algebra
on the complex semisimple Lie group $G$ with respect
to the so-called standard Poisson bracket.
This fact invokes a further topic: equivariant quantization
of actions of complex semisimple Lie groups as actons of
$U_q(\g)$, or equivalently coactions of $\cl{O}_q(G)$, on quantized spaces. One of important researches on this direction is
Podle\'{s}'s study \cite{MR0919322} on equivariant
quantizations of $2$-sphere. Nowadays his work can be
understood as an example of quantized compact
symmetric spaces, initiated by Letzter (\cite{MR1913438}).

Nontheless, in this paper, we would like to explore an alternative approach to generalizing the Podle\'{s} quantum spheres, focusing on quantum semisimple coadjoint orbits, or equivalently quantum flag manifolds.
The latter terminology is used in papers considering the algebra of
smooth or continuous functions on partial flag manifolds with
the action of the compact real form $K$ of $G$, for instance
(\cite{MR3376147, MR1697598}). In this case a $\ast$-structure also seems to be taken in consideration.
On the other hand, the former is used by researchers considering
$L\backslash G$ with a Levi subgroup $L$ and the algebra $\cl{O}(L\backslash G)$ of algebraic functions on $L\backslash G$, usually which is
not equipped with any $\ast$-structure.

Quantization of such actions has already been studied by several
researchers.
In \cite{MR1817512} Donin gives an explicit description of
the moduli of all Poisson structures compatible with the action and also shows that
such deformation quantizaions can be classified by the set of
formal paths on the moduli space.
We explain this fact in more detail. Fix a Cartan subalgebra $\h$
of the Lie algebra $\g$ and consider the associated root system $\roots$
with a positive system $\roots^+$. The set of simple roots is denoted by
$\simples$. We also choose a subset $S$ of $\simples$, from which
we can generate a closed subsystem $\roots_S$ and a Levi subgroup $L_S$.
In this setting, Donin shows that the set of Poisson brackets on $L_S \backslash G$
compatible with the action of the Poisson-Lie group $G$ can be identified
with the set $X_{L_S \backslash G}$ of $\phi = (\phi_{\alpha})_{\alpha \in \roots\setminus\roots_S} \in \C^{\roots\setminus\roots_S}$ satisfying the following relations:
\begin{enumerate}
 \item $\phi_{-\alpha} = -\phi_{\alpha}$ for $\alpha \in \roots\setminus\roots_S$,
 \item $\phi_{\alpha}\phi_{\beta} + 1 = \phi_{\alpha + \beta}(\phi_{\alpha} + \phi_{\beta})$ when $\alpha,\beta,\alpha + \beta \in \roots\setminus\roots_S$,
 \item $\phi_{\alpha} = \phi_{\beta}$ when $\alpha, \beta \in \roots\setminus\roots_S$ and $\alpha - \beta \in \roots_S$.
\end{enumerate}
He also shows the existence of a holomorphic family
$\{\cl{O}_{h,\phi}(L_S \backslash G)\}_{\phi \in X_{L_S \backslash G}}$ of
equivariant deformation quantizations and that
any equivariant deformation quantization is equivalent to $\cl{O}_{h,\phi(h)}(L_S \backslash G)$, where $\phi(h) \in \C\bbh^{\roots\setminus\roots_S}$
satisfies the same relations (i), (ii), (iii). We would like to emphasize here that his result was establised as a
consequence of some vanishing theorems of relevant cohomological
obstructions. In particular any explicit construction of
those deformations was not presented, at least in the paper.

If $\phi \in X_{L_S \backslash G}$ is generic in the sense
that $\phi_{\alpha} \neq -1$ for all $\alpha \in \roots\setminus\roots_S$, Enriquez-Etingof-Marshall (\cite{MR2126485, MR2349621}) and Mudrov (\cite{MR2304470}) give an explicit construction of
corresponding equivariant deformation quantizations. Namely they
construct left $\rep^{\fin}_h \g$-module categories using the parabolic
induction twisted by $\C\bbh_{\lambda/h}$, where $\lambda \in \h^*$
is an appropriate weight. It should be mentioned that this approach
also appears in the investigation of equivariant quantization with respect
to the trivial Poisson-Lie structure (\cite{MR1952112,MR2182701,MR2141466}).

De Commer (\cite{MR3208147}) also followed these pieces of research, with a slight and significant modification, which enables us to discuss quantizations with respect to $\phi$ with $\phi_{\alpha} \neq -1$ for $\alpha \in \roots^+\setminus\roots_S^+$. There he uses a deformed quantum
enveloping algebra instead of the parabolic induction twisted by $\C\bbh_{\lambda/h}$,
which cannot be formulated in this case as the generic case.
He also discusses some operator algebraic aspects, motivated by construction of real semisimple quantum groups as locally compact quantum groups.

\subsection*{Contents of this paper}
The main purpose is to give an explicit way to construct an
algebraic family of $2$-cocyle twists of left $\rep_h^{\fin}\g$-module categories $\rep_h^{\fin}\frl_S$, parametrized by $X_{L_S \backslash G}$. To achieve
this we define a deformed quantum enveloping algebra (Definition \ref{defn:dqea}), which is
inspired by De Commer's work. Namely such a deformed quantum enveloping
algebra is defined whenever we fix a positive system $\roots_0^+$ and
a character on $2Q_0^-$, the cone generated by $-2\roots_0^+$.
If $\roots_0^+ = \roots^+$, this definition corresponds to De Commer's definition.
It is worth to remark here that the PBW theorem still holds for
such algebras, which involves some observations on quantum PBW vectors
in the usual quantum enveloping algebras. This enables us to
investigate the twisted parabolic induction defined in Subsection 3.3.
Any construction in Section 3 is conducted in the integral setting of
Lusztig's sense.

In Section 4, we discuss a kind of semisimplicity of
the twisted parabolic induction in the formal setting.
Using this fact we construct left $\rep_h^{\fin} \g$-module categories
which are $2$-cocycle twists of $\rep_h^{\fin} \frl_S$. We also
construct quantum semisimple coadjoint orbits using these $2$-cocycles and
determine their semi-classical limits.

In Section 5, we discuss the compatibility of the $2$-cocycles constructed
from different choices of $\roots_0^+$.
At the beginning we discuss an immersion of the moduli space $X_{L_S\backslash G}$ to
the toric scheme $X_{\roots}$ associated to the root system $\roots^+$.
Then we construct a sheaf $\rsfs{U}_{q,X_{\roots}}(\gt)$
of deformed quantum enveloping algebras on $X_{\roots}$,
which also induces a sheaf $\rsfs{U}_{q,X_{L_S\backslash G}}(\gt)$.
This fact immediately implies that the $2$-cocycles constructed in Section $4$
can be glued up to a $2$-cocycle with coefficients in $\cl{O}(X_{L_S\backslash G})$.

Finally, in Section 6, we compare quantum semisimple coadjoint orbits
constructed in Section 4 and those induced from quantum semisimple coadjoint orbits of
intermediate Levi subgroups. Namely we show that
an induction functor becomes a left $\rep_h^{\fin}\g$-module functor
in the specific situation.




\section{Preliminaries}
In this paper $\g$ and $\h$ denote a complex semisimple Lie algebra and its Cartan subalgebra respectively.
The associated set of roots is denoted by $\roots$, which forms a root system with a bilinear form $(\tend,\tend)$ on $\h^*$ induced by a normalized Killing form so that $(\alpha,\alpha) = 2$ for a short root $\alpha$.
The reflection with respect to $\alpha \in \roots$ is
denoted by $s_{\alpha}$. The associated Weyl group is denoted by $W$.

We fix a positive system $\roots^+$, which induces
a triangular decomposition $\g = \n^-\oplus \h\oplus \n^+$
and defines a set $\simples = \{\epsilon_1,\epsilon_2,\dots,\epsilon_r\}$ of simple roots. Here $r$, the number of simple roots, is
the rank of $\g$. We also set $N = \abs{\roots^+}$.
The set of reflections with respect to simple roots generates $W$,
and defines the length function $\map{\ell}{W}{\Z_{\ge 0}}$.
The unique longest element is denoted by $w_0$, whose length is $N$.

We set $d_{\alpha}, \alpha^{\vee}, a_{ij}$ as follows:
\[
d_{\alpha} = \frac{(\alpha,\alpha)}{2},\quad
\alpha^{\vee} = d_{\alpha}^{-1}\alpha,\quad
a_{ij} = (\delta_i^{\vee},\delta_j).
\]
The fundamental weights, which are dual to $(\epsilon_i^{\vee})_i$
with respect to $(\tend,\tend)$, are denoted by $\varpi_i$.
The root lattice $Q$ (resp. $P$) is $\Z$-linear span of $\simples$
(resp. $(\varpi_i)_i$). We also use the positive cone $Q^+$ and $P^+$:
\begin{align*}
 Q^+ &= \Z_{\ge 0}\epsilon_1 + \Z_{\ge 0}\epsilon_2 + \cdots + \Z_{\ge 0}\epsilon_r,\\
 P^+ &= \Z_{\ge 0}\varpi_1 + \Z_{\ge 0}\varpi_2 + \cdots + \Z_{\ge 0}\varpi_r.
\end{align*}
We usually replace $\epsilon_i$ by the symbol $i$ when $\epsilon_i$
appears as a subscript. For instance
we use $s_i, d_i, K_i$ instead of using $s_{\epsilon_i}, d_{\epsilon_i}, K_{\epsilon_i}$.

For $S \subset \simples$, we define $\roots_S$ as $\roots\cap \aspan_{\Z}S$ and $\roots_S^+$ as $\roots^+\cap \roots_S$. The number of positive roots in $\roots_S$ is denote by $N_S$.
Then we have Lie subalgebras $\n_S^{\pm} = \bigoplus_{\alpha \in \roots_S^+} \g_{\pm\alpha}$, $\frb_S^{\pm} = \h\oplus\n_S^{\pm}$, $\frl_S = \n_S^-\oplus\h\oplus\n_S^+$, 
$\fru_S^{\pm} = \bigoplus_{\alpha \in \roots^+\setminus\roots_S^+} \g_{\pm\alpha}$, $\frp_S = \frl_S\oplus \fru_S^+$.
We also define $Q_S, P_S, Q_S^+, P_S^+$ for each subset $S \subset \simples$.

For $q$-integers, we use the following symbols:
\[
 q_{\alpha} = q^{d_{\alpha}},\quad
 [n]_q = \frac{q^n - q^{-n}}{q - q^{-1}},\quad
 [n]_q! = [1]_q[2]_q\cdots [n]_q,\quad 
\qbin{n}{k}{q} = \frac{[n]_q!}{[k]_q![n - k]_q!}
\]
A quantum commutator is defined for vectors in suitable algebras which
admits a weight space decomposition as stated in the next section:
\begin{align*}
 [x,y]_q = xy - q^{-(\wt{x},\wt{y})}yx.
\end{align*}
Finally we prepare some notations on a multi-index $\Lambda = (\lambda_i)_i \in \Z_{\ge 0}^n$.

\begin{itemize}
 \item $\abs{\Lambda} = \sum_{i} \lambda_i$.
 \item $\support \Lambda = \{i\mid \lambda_i \neq 0\}$.
 \item $\Lambda \subset (k,l) \deff \support \Lambda \subset \{k + 1, k + 2,\cdots, l - 1\}$. For an interval $I$, like $[k,l]$,
$\Lambda \subset I$ is defined in a similar way.
 \item $\Lambda < k \deff \Lambda \subset (0,k)$. Similarly
$\Lambda \le k, \Lambda > k, \Lambda \ge k$ are defined.
 \item $\Lambda\cdot\alpha = \sum_i \lambda_i\alpha_i$ for
a sequence $(\alpha_i)_i$ of vectors.
 \item $x^{\Lambda} = x_1^{\lambda_1}x_2^{\lambda_2}\cdots x_n^{\lambda_n}$ for a sequence $(x_i)_i$ in a (possibly non-commutative) ring.
\end{itemize}

\subsection{Quantum enveloping algebras of semisimple Lie algebras}
\label{subsec:quantum group}
Basically we refer the convension in \cite{MR4162277} and \cite{MR1492989}.
Though our main purpose is to construct deformation quantizations,
we will present a definition of the deformed quantum enveloping algebra
in the integral form for generality. Therefore
begin with the definition of the Drinfeld-Jimbo deformation over
$\Q(s)$, the function field of one variables over $\Q$
with a deformation parameter $q = s^{2L}$ where $L$ is
the smallest positive integer such that $(\lambda,\mu) \in L^{-1}\Z$
 for any $\lambda,\mu \in P$.
We set $q^{r} := s^{2Lr}$ for $r \in (2L)^{-1}\Z$ and $q_i = q^{d_i}$.
The \emph{Drinfeld-Jimbo deformation} of $\g$ is a Hopf algebra
$U_q(\g)$ generated by $E_i, F_i, K_{\lambda}$ for $1 \le i \le r$ 
and $\lambda \in P$, with relations
\begin{align*}
 &K_0 = 1, &&K_{\lambda}E_iK_{\lambda}^{-1} = q^{(\lambda,\epsilon_i)}E_i, 
&& [E_i, F_j] = \delta_{ij}\frac{K_i - K_i^{-1}}{q_i - q_i^{-1}},\\
 &K_{\lambda}K_{\mu} = K_{\lambda + \mu}, 
&&K_{\lambda}F_iK_{\lambda}^{-1} = q^{-(\lambda,\epsilon_i)}F_i,
\end{align*}
and the quantum Serre relations:
\begin{align*}
&\sum_{k = 0}^{1 - a_{ij}}(-1)^k\qbin{1 - a_{ij}}{k}{q_i}E_i^{1 - a_{ij} - k}E_jE_i^k = 0,\\
&\sum_{k = 0}^{1 - a_{ij}}(-1)^k\qbin{1 - a_{ij}}{k}{q_i}F_i^{1 - a_{ij} - k}F_jF_i^k = 0.
\end{align*}
The coproduct $\Delta$, the antipode $S$ and the counit $\eps$ 
are given as follows on the generators:
\begin{align*}
 &\Delta(K_{\lambda}) = K_{\lambda}\tensor K_{\lambda}, &&S(K_{\lambda}) = K_{\lambda}^{-1}, &&\eps(K_{\lambda}) = 1,\\
 &\Delta(E_i) = E_i\tensor K_i + 1\tensor E_i, &&S(E_i) = -E_iK_i^{-1}, &&\eps(E_i) = 0,\\
 &\Delta(F_i) = F_i\tensor 1 + K_i^{-1}\tensor F_i, &&S(F_i) = -K_iF_i, &&\eps(F_i) = 0.
\end{align*}


Next we introduce some subalgebras of $U_q(\g)$. 
The most fundamental ones are $U_q(\n^+), U_q(\n^-)$ and 
$U_q(\h)$, which are genereted by 
$E_i, F_i, K_{\lambda}$ respectively. These allow us to
decompose $U_q(\g)$ into the tensor products 
$U_q(\n^{\pm})\tensor U_q(\h)\tensor U_q(\n^{\mp})$ via the 
multiplication maps. We also use $U_q(\frb^{\pm})$ for the 
subalgebras generated by $U_q(\h)$ and $U_q(\n^{\pm})$ respectively.
Note that $U_q(\frb^{\pm})$ are Hopf subalgebras of $U_q(\g)$.

For a $U_q(\h)$-module $M$ and $v \in M$, we say that
$v$ is a weight vector of weight $\lambda \in P$ when
$K_{\mu}v = q^{(\mu,\lambda)}v$ for all $\mu \in P$. In this case
$\lambda$ is denoted by $\wt v$. The submodule of elements of
weight $\lambda$ is denoted by $M_{\lambda}$.
To consider the weight of
an element of $U_q(\g)$, we regard $U_q(\g)$ as a
$U_q(\h)$-module by the restriction of the left adjoint action
$x\triangleright y = x_{(1)}yS(x_{(2)})$.

Next we describe the braid group action on $U_q(\g)$ and
the quantum PBW bases.
At first we have an algebra automorphism
$\cl{T}_i$ on $U_q(\g)$ for each $\delta_i \in \simples$,
which satisfies
\begin{align*}
\cl{T}_i(K_{\lambda}) = K_{s_i(\lambda)}, \quad
\cl{T}_i(E_i) = -K_iF_i, \quad
\cl{T}_i(F_i) = -E_iK_i^{-1}
\end{align*}
and other formulae in \cite[Theorem 3.58]{MR4162277}
which determine $\cl{T}_i$ uniquely.

Then the family $(\cl{T}_i)_i$ satisfies the Coxeter relations and
defines a braid group action on $U_q(\g)$. Especially we have
$\cl{T}_w$ for each $w \in W$, which is given by
$\cl{T}_w = \cl{T}_{i_1}\cl{T}_{i_2}\cdots\cl{T}_{i_{\ell(w)}}$ where
$w = s_{i_1}s_{i_2}\cdots s_{i_{\ell(w)}}$ is a reduced expression.

This action produces PBW bases of $U_q(\g)$.
Let $w_0$ be the longest element in $W$ and fix its reduced expression
$w_0 = s_{\bs{i}} = s_{i_1}s_{i_2}\cdots s_{i_N}$, where
$\bs{i} = (i_1,i_2,\dots,i_N)$.
Then each $\alpha \in R^+$ has a unique positive integer $k \le N$
with $\alpha = \alpha^{\bs{i}}_{k} := s_{i_1}s_{i_2}\cdots s_{i_{k - 1}}(\eps_{i_k})$.
Finally we set $E_{\bs{i},\alpha}$ and 
$F_{\bs{i},\alpha}$, the quantum root vectors, as follows:
\begin{align*}
 E_{\bs{i},\alpha} &= E_{\bs{i},k} 
:= \cl{T}_{s_{i_1}s_{i_2}\cdots s_{i_{k - 1}}}(E_{i_k}) 
= \cl{T}_{s_{i_1}}\cl{T}_{s_{i_2}}\cdots\cl{T}_{s_{i_{k - 1}}}(E_{i_k}), \\
 F_{\bs{i},\alpha} &= F_{\bs{i},k} 
:= \cl{T}_{s_{i_1}s_{i_2}\cdots s_{i_{k - 1}}}(F_{i_k}) 
= \cl{T}_{s_{i_1}}\cl{T}_{s_{i_2}}\cdots\cl{T}_{s_{i_{k - 1}}}(F_{i_k}).
\end{align*}
Though these elements depend on $\bs{i}$, we still have an
analogue of the Poincar\'{e}-Birkhoff-Witt theorem in $U_q(\g)$ 
i.e. $\{F_{\bs{i}}^{\Lambda^-}K_{\mu}E_{\bs{i}}^{\Lambda^+}\}_{\Lambda^{\pm},\mu}$ forms a basis of $U_q(\g)$. Each element of
this family is called a \emph{quantum PBW vector}.

\subsection{The Drinfeld pairing and the quantum Killing form}
Next we recall the Drinfeld pairing and the quantum Killing form. 
The \emph{Drinfeld pairing} is a unique bilinear form
$\map{\tau}{U_q(\frb^+)\times U_q(\frb^-)}{k}$ with
\begin{align*}
 \tau(K_{\lambda}, K_{\mu}) = q^{-(\lambda,\mu)},\quad
 \tau(E_i, K_{\lambda}) = 0 = \tau(F_i, K_{\mu}),\quad
 \tau(E_i, F_j) = -\frac{\delta_{ij}}{q_i - q_i^{-1}}
\end{align*}
and the following conditions:
\begin{align*}
 &\tau(xx',y) = \tau(x,y_{(1)})\tau(x',y_{(2)}),\quad
 &&\tau(x,1) = \eps(x),\quad \\
 &\tau(x,yy') = \tau(x_{(1)},y')\tau(x_{(2)},y),\quad
 &&\tau(1,y) = \eps(y).
\end{align*}
Here $\Delta(x) = x_{(1)}\tensor x_{(2)}$ is Sweedler's notation.

The PBW bases are orthogonal with respect to this form:
\begin{align*}
 \tau(E_{\bs{i}}^{\Lambda^+}, F_{\bs{i}}^{\Lambda^-})
= \delta_{\Lambda^+ \Lambda^-}\prod_{k = 1}^N
(-1)^{\lambda_k^+}
q_{\alpha^{\bs{i}}_k}^{\lambda^+_k(\lambda^+_k - 1)/2}
\frac{[\lambda^+_k]_{q_{\alpha^{\bs{i}}_k}!}}{(q_{\alpha^{\bs{i}}_k} - q^{-1}_{\alpha^{\bs{i}}_k})^{\lambda^+_k}}.
\end{align*}
The following formulae will be used frequently:
\begin{align}\label{eq:translation by K}
&\tau(xK_{\mu},yK_{\nu}) = q^{-(\mu,\nu)}\tau(x,y),
&&x \in U_q(\n^+), y \in U_q(\n^-), \\ 
\label{eq:product via the Killing form}
&xy = \tau(S(x_{(1)}),y_{(1)})y_{(2)}x_{(2)}\tau(x_{(3)},y_{(3)}),
&&x \in U_q(\frb^+), y \in U_q(\frb^-).
\end{align}

The \emph{quantum Killing form}
$\map{\kappa}{U_q(\g)\times U_q(\g)}{k}$ is given 
by the following formula:
\begin{align*}
 \kappa(X^-K_{\mu}S^{-1}(X^+),Y^-K_{\nu}S(Y^+)) = q^{(\mu,\nu)/2}\tau(X^+,Y^-)\tau(Y^+,X^-).
\end{align*}
Ad-invariance, an important property of the classical
Killing form, still holds in this case:
\[
 \kappa(Z\triangleright X,Y) = \kappa(X,S(Z)\triangleright Y)
\]

\subsection{The restricted integral form}

At the end of this section we recall the restricted integral form $U_q^{\A}(\g)$.
Set $\A := \Z[s^2,s^{-2}] \subset \Q(s)$.
Then $U_q(\g)$ has a Hopf subalgebra $U_q^{\A}(\g)$ over $\A$,
called the restricted integral form,
which is generated by the following elements:
\begin{align*}
 E_i^{(r)} = \frac{1}{[r]_{q_i}!}E_i^r,\quad
 F_i^{(r)} = \frac{1}{[r]_{q_i}!}F_i^r,\quad
 K_{\lambda},\quad
 [K_i;0]_{q_i} = \frac{K_i - K_i^{-1}}{q_i - q_i^{-1}}.
\end{align*}
The upper (resp. lower) triangular subalgebra $U_q^{\A}(\n^{\pm})$
are defined similarly to the non-integral case.
They are free over $\A$ since the following elements 
form a basis:
\begin{align*}
 E_{\bs{i}}^{(\Lambda)}
= E_{\bs{i},1}^{(\lambda_1)}E_{\bs{i},2}^{(\lambda_2)}\cdots
E_{\bs{i},N}^{(\lambda_N)},\quad
 F_{\bs{i}}^{(\Lambda)} 
= F_{\bs{i},1}^{(\lambda_1)}
F_{\bs{i},2}^{(\lambda_2)}
\cdots
F_{\bs{i},N}^{(\lambda_N)}
\end{align*}
We also have the Cartan subalgebra and the Borel subalgebras, defined
as $U_q^{\A}(\h) = U_q(\h)\cap U_q^{\A}(\g)$ and $U_q^{\A}(\frb^{\pm}) = U_q^{\A}(\n^{\pm})U_q^{\A}(\g)$ respectively.

Finally, for a commutative $\A$-algebra $k$, we define $U_q^k(\g)$
as $k\tensor_{\A} U_q^{\A}(\g)$. we also define
$U_q^{k}(\n^{\pm}), U_q^{k}(\h)$ and $U_q^{k}(\frb^{\pm})$
in a similar way.

\section{Deformed quantum enveloping algebras}
\label{sect:deformed QEA}
\subsection{More on quantum PBW bases}
Before defining deformed quantum enveloping algebras,
we collect some technical facts on quantum PBW vectors in the usual 
quantum enveloping algebra.
For a reduced expression $s_{\bs{i}}$ of $w_0$, we set
$\aE_{\bs{i},\alpha} = E_{\bs{i},\alpha}K_{\alpha}^{-1}, \aK_{\mu} = K_{\mu}$ and $\aF_{\bs{i},\alpha} = F_{\bs{i},\alpha}$.
\begin{rema}
In the following, various expansion formulae
with respect to the modified quantum
PBW vectors $\aF_{\bs{i}}^{(\Lambda^-)}K_{\mu}\aE_{\bs{i}}^{(\Lambda^+)}$ are stated.
We sometimes apply them to other generators, like $\hE_{\bs{i},\alpha}, \hK_{\mu}, \hF_{\bs{i},\alpha}$.
Moreover we use their variations: for instance,
using $\hE_{\bs{i}}^{\Lambda^+}$ instead of $\hE_{\bs{i}}^{(\Lambda^+)}$. These application can be validated by considering the ``easiest'' case,
in which these generators are scalar multiple of the original generators. Then the other cases are reduced to this case by changing the base ring.
\end{rema}

The following lemma is fundamental in this subsection.

\begin{lemm} \label{lemm:rotation lemma}
Let $s_{\bs{i}}$ be a reduced expression of $w_0$
and set $w = s_{i_1}s_{i_2}\cdots s_{i_k}$ for a fixed $k$
with $1 \le k \le N$. Then there is another reduced expression 
$s_{\bs{j}}$ which satisfies the following:
\begin{align*}
 \cl{T}_w(\aE_{\bs{j},l}) &= 
\begin{cases}
 \aE_{\bs{i},k + l} & (1 \le l \le N - k) \\
 -q_{\alpha_{l - (N - k)}^{\bs{i}}}^{-2}\aF_{\bs{i},l - (N - k)}\aK_{2\alpha^{\bs{i}}_{l - (N - k)}} & (N - k < l \le N),
\end{cases}
\\
 \cl{T}_w(\aF_{\bs{j},l}) &= 
\begin{cases}
 \aF_{\bs{i},k + l} & (1 \le l \le N - k) \\
 -\aE_{\bs{i},l - (N - k)} & (N - k < l \le N),
\end{cases}
\end{align*}
\end{lemm}
\begin{proof}
At first we assume $k = 1$. Set
$\epsilon = s_{i_N}s_{i_{N - 1}}\cdots s_{i_2}(\epsilon_{i_1})$.
Then this is a simple root and
$s_{\bs{j}} = s_{i_2}s_{i_3}\cdots s_{i_N}s_{\eps}$
is a reduced expression of $w_0$, which has the required property.
Actually we can see from the construction that
$\cl{T}_{i_1}(\aE_{\bs{j},l}) = \aE_{\bs{i},l + 1}$ and
$\cl{T}_{i_1}(\aF_{\bs{j},l}) = \aF_{\bs{i},l + 1}$ for
$1 \le l \le N-1$.

For the case $l = N$,
$\aE_{\bs{j},N} = \aE_{i_1}$ and $\aF_{\bs{j},N} = \aF_{i_1}$ hold
since $\alpha_N^{\bs{j}} = \epsilon_{i_1}$ is a simple root.
Hence we have
\begin{align*}
 \cl{T}_{i_1}(\aE_{\bs{j},N})
&= \cl{T}_{i_1}(\aE_{i_1}) 
= -K_{i_1}\aF_{i_1}K_{i_1}
= -q_{\alpha^{\bs{i}}_1}^{-2}\aF_{\bs{i},1}\aK_{2\alpha^{\bs{i}}_1},\\
 \cl{T}_{i_1}(\aF_{\bs{j},N}) 
&= \cl{T}_{i_1}(\aF_{i_1}) 
= -\aE_{i_1}
= -\aE_{\bs{i},1}.
\end{align*}

Now we show the statement by induction on $k$. Assume 
the statement with $k = k_0$ holds for all reduced expressions. Then
we can take $\bs{i'}$ for $\bs{i}$ and $k = k_0$, which automatically
satisfies $i'_1 = i_{k + 1}$. Then we can take $\bs{j}$ for $\bs{i'}$
and $k = 1$, which is the required one for $\bs{i}$ and $k = k_0 + 1$.
\end{proof}

Next we investigate quantum commutators, whose definition
in this paper is $[x,y]_q := xy - q^{-(\wt x, \wt y)}yx$
for weight vectors $x,y$.

The following is an easy consequence of a well-known result,
called the \emph{Levend\"{o}rskii-Soibelman relation} in literature.

\begin{prop}[c.f. {\cite[Proposition 5.5.2]{MR1116413}}]
\label{prop:LS relation}
For $1 \le k < l \le N$, we have the following expansion with $C_{\Lambda}^{\pm} \in \A$:
\begin{align*}
 [\aE_{\bs{i},l},\aE_{\bs{i},k}]_q
= \sum_{\substack{\Lambda \subset (k,l)\\\Lambda\cdot\alpha^{\bs{i}} = \alpha^{\bs{i}}_k + \alpha^{\bs{i}}_l}}
 C_{\Lambda}^+\aE_{\bs{i}}^{(\Lambda)}, \quad
 [\aF_{\bs{i},l},\aF_{\bs{i},k}]_q
= \sum_{\substack{\Lambda \subset (k,l)\\\Lambda\cdot\alpha^{\bs{i}} = \alpha^{\bs{i}}_k + \alpha^{\bs{i}}_l}}
 C_{\Lambda}^-\aF_{\bs{i}}^{(\Lambda)},
\end{align*}
\end{prop}

Combining this with Lemma \ref{lemm:rotation lemma},
we can obtain a ``mixed'' version of these formulae.

\begin{prop} \label{prop:mixed LS relation}
Fix a reduced expression $s_{\bs{i}}$ of the longest element $w_0$.
and positive integers $1 \le k < l \le N$.
\begin{enumerate}
 \item The following $\A$-submodules are closed under multiplication:
\begin{align*}
&\aspan_{\A}\{\aF_{\bs{i}}^{\Lambda^-}\aE_{\bs{i}}^{\Lambda^+}\mid \Lambda^+ \le k < l \le \Lambda^-\}, \\
&\aspan_{\A}\{\aE_{\bs{i}}^{\Lambda^+}\aK_{\mu}\aF_{\bs{i}}^{\Lambda^-}\mid \mu \in P,\,\Lambda^- \le l < k \le \Lambda^+\}.
\end{align*}
 \item The latter subalgebra coincides with the following:
\begin{align*}
 \aspan_{\A}\{\aF_{\bs{i}}^{(\Lambda^-)}\aK_{\mu}\aE_{\bs{i}}^{(\Lambda^+)}\mid \mu \in P, \Lambda^- \le l < k \le \Lambda^+\}.
\end{align*}
 \item We have the following expansion with $C_{k,l}(\Lambda^-,\Lambda^+), C_{l,k}(\Lambda^-,\Lambda^+) \in \A$:
\begin{align} \label{eq:mixed LS relation}
[\aE_{\bs{i},k},\aF_{\bs{i},l}]_q &=
  \sum_{\substack{\Lambda^+ <k < l < \Lambda^-,\\ \mu :=\alpha_k^{\bs{i}} - \Lambda^+\cdot\alpha^{\bs{i}} = \alpha_l^{\bs{i}} - \Lambda^-\cdot\alpha^{\bs{i}} \in Q^+}}
C_{k,l}(\Lambda^-,\Lambda^+) \aF_{\bs{i}}^{(\Lambda^-)}\aE_{\bs{i}}^{(\Lambda^+)}, \\
[\aE_{\bs{i},l},\aF_{\bs{i},k}]_q &= \sum_{\substack{\Lambda^- < k < l < \Lambda^+\\ \mu := \alpha_k^{\bs{i}} - \Lambda^+\cdot\alpha^{\bs{i}} = \alpha_l^{\bs{i}} - \Lambda^-\cdot\alpha^{\bs{i}} \in Q^+}} 
C_{l,k}(\Lambda^-,\Lambda^+) \aF_{\bs{i}}^{(\Lambda^-)}\aK_{2\mu}^{-1}\aE_{\bs{i}}^{(\Lambda^+)}.
\end{align}
\end{enumerate}
\end{prop}
\begin{proof}
(i) This follows from Lemma \ref{lemm:rotation lemma} immediately.

(ii) Note that the latter subalgebra in (i) is generated by $(\aE_{\bs{i},n})_{n \le l}$, $(\aF_{\bs{i},m})_{m \ge k}$ and $(\aK_{\mu})_{\mu \in P}$. Hence it suffices to show that the following $\A$-submodules are closed under multiplication for all $0 \le m \le l$:
\begin{align*}
\cl{I}_m := \aspan_{\A}\{\aF_{\bs{i}}^{(\Lambda^-_1)}\aE_{\bs{i}}^{(\Lambda^+)}\aK_{\mu}\aF_{\bs{i}}^{(\Lambda^-_2)}\mid \mu \in P, \Lambda^-_2 \le m < \Lambda^-_1 \le l < k \le \Lambda^+\}.
\end{align*}
If $m = l$, there is nothing to prove. For the downward induction step, assume that the statement holds for a fixed $m$.
It suffices to show $\cl{I}_{m - 1}$ is invariant under the right multiplication by $\aE_{\bs{i},m}^{(r)}$ for all $r \ge 1$. There is nothing to prove when $r = 0$. For a general $r$, we can see $\cl{I}_{m - 1}\aE_{\bs{i},m}^{(r)} \subset \cl{I}_{m - 1}\aE_{\bs{i},m}^{(r - 1)}$ using Lemma \ref{lemm:rotation lemma} and
Proposition \ref{prop:LS relation}. Hence an induction argument works.

(iii) Take $\bs{j}$ using Lemma \ref{lemm:rotation lemma} for $\bs{i}$ and $k$. Then we have
\begin{align*}
\aE_{\bs{i},k}\aF_{\bs{i},l} - q^{-(\alpha^{\bs{i}}_k,\alpha^{\bs{i}}_l)}\aF_{\bs{i},l}\aE_{\bs{i},k}
= q^{(\alpha_k^{\bs{i}},\alpha_l^{\bs{i}})}\cl{T}_w(\aF_{\bs{j},l - k}\aF_{\bs{j},n} - q^{(\alpha_{l - k}^{\bs{j}},\alpha_n^{\bs{j}})}\aF_{\bs{j},n}\aF_{\bs{j},l-k}).
\end{align*}
Now Proposition \ref{prop:LS relation} leads us to the desired expansion without the constraint on $\mu$.
A similar argument works for the other case.

Finally we show the positivity of $\mu$. Since
$\Delta(F_i) = F_i \tensor 1 + K_i^{-1}\tensor F_i$ holds,
$\Delta(F_{\bs{i},\alpha})$ is an $\A$-linear combination of
$F_{\bs{i}}^{(\Lambda_1)}K_{\Lambda_2\cdot \alpha^{\bs{i}}}^{-1}\tensor F_{\bs{i}}^{(\Lambda_2)}$ with $\Lambda_1\cdot\alpha^{\bs{i}} + \Lambda_2\cdot\alpha^{\bs{i}} = \alpha$.
Using this fact twice we can see that
$(\Delta\tensor\id)\Delta(\aF_{\bs{i},\alpha})$ is an $\A$-linear
combination of $F_{\bs{i}}^{(\Lambda_1)}K_{\Lambda_2\cdot\alpha^{\bs{i}}+ \Lambda_3\cdot\alpha^{\bs{i}}}^{-1}\tensor \aF_{\bs{i}}^{(\Lambda_2)}\aK_{\Lambda_3\cdot\alpha^{\bs{i}}}^{-1}\tensor F_{\bs{i}}^{(\Lambda_3)}$
with $\Lambda_1\cdot\alpha^{\bs{i}} + \Lambda_2\cdot\alpha^{\bs{i}} + \Lambda_3\cdot\alpha^{\bs{i}} = \alpha$.
A similar fact holds for $(\Delta\tensor\id)\Delta(\aE_{\bs{i},\beta})$.
Then the formula (\ref{eq:product via the Killing form}) implies that
$\aE_{\bs{i},\alpha^+}\aF_{\bs{i},\alpha^-}$ is an $\A$-linear combination of elements of the following form:
\begin{align*}
\tau(S(E_{\bs{i}}^{(\Lambda_1^+)}K_{\alpha^+}^{-1}), F_{\bs{i}}^{(\Lambda_1^-)}K_{\Lambda_2^-\cdot\alpha^{\bs{i}} + \Lambda_3^-\cdot\alpha^{\bs{i}}}^{-1})\tau(E_{\bs{i}}^{(\Lambda_3^+)}K_{\Lambda_3^+\cdot\alpha^{\bs{i}}}^{-1},F_{\bs{i}}^{(\Lambda_3^-)})
\aF_{\bs{i}}^{(\Lambda_2^-)}\aK_{\Lambda_3^-\cdot\alpha^{\bs{i}} + \Lambda_3^+\cdot\alpha^{\bs{i}}}^{-1}\aE_{\bs{i}}^{(\Lambda_2^+)}
\end{align*}
with $\Lambda_1^{\pm}\cdot\alpha^{\bs{i}} + \Lambda_2^{\pm}\cdot\alpha^{\bs{i}} + \Lambda_3^{\pm}\cdot\alpha^{\bs{i}} = \alpha^{\pm}$.
Since this term does not vanish only when $\Lambda_3^+ = \Lambda_3^-$,
we have $\mu = \Lambda_3^{\pm}\cdot\alpha^{\bs{i}} \in Q^+$.
\end{proof}

As an application we can obtain formulae on
the coproduct.

\begin{prop} \label{prop:expansion of coproducts}
Fix a reduced expression $s_{\bs{i}}$ of $w_0$.
Then we have the following expansion with $C_k^{\pm}(\Lambda^l,\Lambda^r) \in \A$:
\begin{align}
\label{eq:coproduct of E}
&\begin{aligned}
 \Delta(\aE_{\bs{i},k}) = E_{\bs{i},k}&K_{\alpha^{\bs{i}}_k}^{-1}\tensor 1 + K_{\alpha^{\bs{i}}_k}^{-1}\tensor \aE_{\bs{i},k}\\
 &+ \sum_{\Lambda_r < k < \Lambda_l, \Lambda_l\cdot\alpha^{\bs{i}} + \Lambda_r\cdot\alpha^{\bs{i}} = \alpha^{\bs{i}}_k} C_k^+(\Lambda_l,\Lambda_r)K_{\alpha^{\bs{i}}_k}^{-1}E_{\bs{i}}^{(\Lambda_l)}\tensor \aE_{\bs{i}}^{(\Lambda_r)},
\end{aligned} \\
\label{eq:coproduct of F}
&\begin{aligned}
 \Delta(\aF_{\bs{i},k}) = F_{\bs{i},k}&\tensor 1 + K_{\alpha^{\bs{i}}_k}^{-1}\tensor \aF_{\bs{i},k}\\
 &+ \sum_{\Lambda_l < k < \Lambda_r, \Lambda_l\cdot\alpha^{\bs{i}} + \Lambda_r\cdot\alpha^{\bs{i}} = \alpha^{\bs{i}}_k} C_k^-(\Lambda_l,\Lambda^r)F_{\bs{i}}^{(\Lambda^l)}K_{\Lambda^r\cdot\alpha^{\bs{i}}}^{-1}\tensor \aF_{\bs{i}}^{(\Lambda^r)}.
\end{aligned}
\end{align}
\end{prop}
We need some preparations to prove this proposition.
Recall the quantum exponential function:
\begin{align*}
 \exp_q x = \sum_{k = 0}^{\infty} \frac{q^{k(k-1)/2}}{[k]_q!} x^k
\end{align*}
Let $s_{\bs{i}}$ be a reduced expression of an element of $W$.
\begin{align*}
\Delta(\cl{T}_i(z))
= A_i(\cl{T}_i\tensor \cl{T}_i)(\Delta(z))A_i^{-1},
\end{align*}
where $A_i := \exp_{q_i} ((q_i - q_i^{-1})K_iF_i\tensor \aE_i)$,
whose inverse is given by
$A_i^{-1} = \exp_{q_i^{-1}} (-(q_i - q_i^{-1})K_iF_i\tensor \aE_i)$.
Applying this formula iteratively to a reduced expression
$s_{\bs{i}}$ of an element $w \in W$, we also have
\begin{align*}
 \Delta(\cl{T}_w(z)) = A_{\bs{i},1}A_{\bs{i},2}\cdots A_{\bs{i},\ell(w)}(\cl{T}_w\tensor \cl{T}_w)(\Delta(z))A_{\bs{i},\ell(w)}^{-1}A_{\bs{i},\ell(w) - 1}^{-1}\cdots A_{\bs{i},1}^{-1}
\end{align*}
where
\begin{align*}
 A_{\bs{i},l} &= \exp_{q_{\alpha^{\bs{i}}_l}}((q_{\alpha^{\bs{i}}_l} - q_{\alpha^{\bs{i}}_l}^{-1})K_{\alpha^{\bs{i}}_l}F_{\bs{i},l}\tensor \aE_{\bs{i},l}), \\
 A_{\bs{i},l}^{-1} &= \exp_{q_{\alpha^{\bs{i}}_l}^{-1}}(-(q_{\alpha^{\bs{i}}_l} - q_{\alpha^{\bs{i}}_l}^{-1})K_{\alpha^{\bs{i}}_l}F_{\bs{i},l}\tensor \aE_{\bs{i},l}).
\end{align*}
\begin{proof}
Let $s_{\bs{i}}$ be a reduced expression of the longest element
compatible with $\roots_0^+$
and $w_k = s_{i_1}s_{i_2}\cdots s_{i_k}$.
Then, for any $1 \le k \le n$, we have
\begin{align*}
 \Delta(\aE_{\bs{i},k})
= A_{\bs{i},1}A_{\bs{i},2}\cdots A_{\bs{i},k-1}(\cl{T}_{w_{k - 1}}\tensor \cl{T}_{w_{k - 1}})(\Delta(\aE_{i_k}))A_{\bs{i},k - 1}^{-1}A_{\bs{i},k - 2}^{-1}\cdots A_{\bs{i},1}^{-1}.
\end{align*}
On the other hand we have
\begin{align*}
 (\cl{T}_{w_{k - 1}}\tensor \cl{T}_{w_{k - 1}})(\Delta(\aE_{i_k}))
= E_{\bs{i},k}K_{\alpha^{\bs{i}}_k}^{-1}\tensor 1 + K_{\alpha^{\bs{i}}_k}^{-1}\tensor \aE_{\bs{i},k}
\end{align*}
Combining with Proposition \ref{prop:LS relation}, Proposition \ref{prop:mixed LS relation} (i) and $\Delta(\aE_{\bs{i},k}) \in U_q^{\A}(\frb^+)\tensor U_q^{\A}(\frb^+)$, we can see (\ref{eq:coproduct of E}).

The same argument works for (\ref{eq:coproduct of F})
since the following formula holds:
\begin{align*}
(\cl{T}_{w_{k - 1}}\tensor \cl{T}_{w_{k - 1}})(\Delta(\aF_{i_k}^{(r)}))
=\sum_{j = 0}^r q_{\alpha^{\bs{i}}_k}^{j(r - j)}F_{\bs{i},k}^{(r - j)}K_{\alpha^{\bs{i}}_k}^{-j}\tensor \aF_{\bs{i},k}^{(j)}.
\end{align*}
For a reference, see \cite[Proof of Lemma 3.18]{MR4162277}.
\end{proof}

\subsection{Deformed quantum enveloping algebra}

We fix a positive system $\roots_0^+$,
possibly different from $\roots^+$.
Then there is a unique element $w \in W$ such that $w(\roots^+) = \roots_0^+$. Note $\roots^+\setminus\roots_0^+ = \{\alpha \in\roots^+ \mid w^{-1}(\alpha) \in -\roots^+\}$ for such $w$.

We say that a reduced expression $s_{\bs{i}}$ of the longest
element is \emph{compatible with $\roots_0^+$} if it begins with
a reduced expression of $w$ i.e. $w = s_{i_1}s_{i_2}\cdots s_{i_{\ell(w)}}$.
Such an expression always exists since $\ell(w_0) = \ell(w) + \ell(w^{-1}w_0)$.
Note that $s_{\bs{i}}$ is compatible with $\roots_0^+$ if and only if
$\roots^+\setminus\roots_0^+ = \{\alpha^{\bs{i}}_{n}\}_{n = 1}^{\ell(w)}$.

For a monoid $M$, the monoid algebra with coefficients in a commutative
ring $k$ is denoted by $k[M]$. It has a canonical $k$-basis
$\{e_{m}\}_{m \in M}$, for which we have $e_me_{m'} = e_{mm'}$.

Let $Q_0^-$ be the additive submonoid of $Q$
generated by $-\roots_0^+$. By the universal property of
monoid algebras, we can identify
an $\A[2Q_0^-]$-algebra with a pair of an $\A$-algebra $k$ and
a monoid homomorphism $\map{\chi}{2Q_0^-}{k}; \lambda \longmapsto \chi_{\lambda}$ with respect to the multiplication on $k$. In light of this fact, we refer such a pair as a commutative $\A[2Q_0^-]$-algebra.
\begin{defn} \label{defn:dqea}
Let $(k,\chi)$ be a commutative $\A[2Q_0^-]$-algebra and
fix a reduced expression $s_{\bs{i}}$ compatible with $\roots_0^+$.
We define $U_{q,\chi}^{k}(\gt)$ as
$k\tensor_{\A[2Q_0^-]} U_{q,e}^{\A[2Q_0^-]}(\gt)$,
where $U_{q,e}^{\A[2Q_0^-]}(\gt)$ is an $\A[2Q_0^-]$-subalgebra of
$U_q^{\A[P]}(\g)$ generated by the following elements:
\begin{itemize}
 \item ${\displaystyle \hF_{\bs{i},\alpha}^{(r)} := 
\begin{cases}
 e_{2\alpha}^r \aF_{\bs{i},\alpha}^{(r)} & \text{for }\alpha \in \roots^+\setminus\roots_0^+, r \ge 0\\
 \aF_{\bs{i},\alpha}^{(r)} & \text{for } \alpha \in \roots^+\cap \roots_0^+, r\ge 0,
\end{cases}}$
 \item $\hK_{\lambda} := e_{-\lambda}\aK_{\lambda}$ for $\lambda \in P$.
 \item $\hE_{\bs{i},\alpha} := (q_{\alpha} - q_{\alpha}^{-1})\aE_{\bs{i},\alpha}$ for $\alpha \in \roots^+$.
\end{itemize}
We also define $U_{q,\chi}^k(\tg)$ by using
$\vF_{\bs{i},\alpha} = (q_{\alpha} - q_{\alpha}^{-1})\hF_{\bs{i},\alpha}, \vK_{\lambda} = \hK_{\lambda}, \vE_{\bs{i},\alpha}^{(r)} = \aE_{\bs{i},\alpha}^{(r)}$.
\end{defn}
\begin{rema}
These definitions actually do not depend on the choice of
$s_{\bs{i}}$ since
$\mathrm{span}_{\A}\{F_{\bs{i}}^{(\Lambda)}\}_{\Lambda \le \abs{\roots^+\setminus\roots_0^+}} = U_q^{\A}(\n^-)\cap \cl{T}_w^{-1}(U_q^{\A}(\frb^+))$ and $\mathrm{span}_{\A}\{F_{\bs{i}}^{(\Lambda)}\}_{\Lambda > \abs{\roots^+\setminus\roots_0^+}} = U_q^{\A}(\n^-)\cap \cl{T}_w(U_q^{\A}(\n^-))$ are independent of $s_{\bs{i}}$, where $w$ is the unique element with $w(\roots^+) = \roots_0^+$.

On the other hand,
though we can consider a subalgebra of $U_q^{\A[P]}(\g)$ generated by
the elements above constructed from a non-compatible reduced expression $s_{\bs{i}}$,
it can be different from $U_{q,e}^{\A[2Q_0^-]}(\gt)$. As an example
we consider $\g = \frak{sl}_3$. The simple roots is denoted by $\alpha$ and $\beta$. Then the reduced expressions of the longest element are
$s_{\bs{i}} = s_{\alpha}s_{\beta}s_{\alpha}$ and $s_{\bs{j}} := s_{\beta}s_{\alpha}s_{\beta}$.
Consider $\roots_0^+ := \{\beta, -\alpha, -\alpha -\beta\}$.
Then $s_{\bs{i}}$ is compatible with $\roots_0^+$ and $s_{\bs{j}}$ is not.
In this case we have
\begin{align*}
 e_{2(\alpha + \beta)}F_{\bs{j}, \alpha + \beta}
= e_{2(\alpha + \beta)}(F_{\beta}F_{\alpha} - q^{-1}F_{\alpha}F_{\beta})
= -q^{-1}\hF_{\bs{i},\alpha + \beta} - e_{2\beta}(q - q^{-1})\hF_{\beta}\hF_{\alpha}.
\end{align*}
Since $e_{2\beta}$ is not invertible in $\A[2Q_0^-]$ and the PBW theorem holds for $U_{q,e}^{\A[2Q_0^-]}(\gt)$, as shown in Proposition \ref{prop:PBW theorem for gt}, $e_{2(\alpha + \beta)}F_{\bs{i},\alpha + \beta}$
is not in $U_{q,e}^{\A[2Q_0^-]}(\gt)$.
\end{rema}
\begin{rema}
We cannot remove $\hF_{\bs{i},\alpha}$ with non-simple root $\alpha$
from the set of generators. As an example, we consider $\frak{sl}_3$ and
a positive system $\roots_0^+ := \{\alpha + \beta,\beta, -\alpha\}$. Then $s_{\bs{i}}$ is compatible with $\roots_0^+$.
Then we can calculate some commutation relations of the generators as follows:
\begin{align*}
 [\hF_{\alpha},\hF_{\beta}]_q = e_{2\alpha}\hF_{\alpha + \beta},\,
 [\hE_{\alpha + \beta},\hF_{\alpha}]_q = \hK_{2\alpha}^{-1}\hE_{\beta},\,
 [\hE_{\alpha + \beta},\hF_{\beta}]_q = -\hE_{\alpha},
\end{align*}
here we omit the subscript $\bs{i}$.
Since $e_{2\alpha}$ is not invertible in $\A[2Q_0^-]$,
the subalgebra generated by $\hF_{\alpha}^{(r)},\hF_{\beta}^{(r)}, \hK_{\lambda}, \hE_{\alpha}$, $\hE_{\alpha + \beta}, \hE_{\beta}$ does not
contain $\hF_{\alpha + \beta}$.
\end{rema}

\begin{prop} \label{prop:well-definedness of coaction}
The canonical left coaction of $U_q^{\A[2Q_0^-]}(\g)$ on
itself restricts to a left coaction on $U_{q,e}^{\A[2Q_0^-]}(\gt)$.
In particular $U_{q,\chi}^{k}(\gt)$ has a canonical left $U_q^k(\g)$-coaction. Similarly $U_{q,\chi}^k(\tg)$ has a canonical left $U_q^k(\g)$-coaction.
\end{prop}

\begin{proof}
Proof of Proposition \ref{prop:expansion of coproducts} shows
the statement since $\Delta(\hF_{\bs{i},\alpha}^{(r)})$,
$\Delta(\hK_{\mu}), \Delta(\hE_{\bs{i},\alpha})$ and
$\pm K_{\alpha}F_{\bs{i},\alpha}\tensor \hE_{\bs{i},\alpha}$,
whose quantum exponential is $A_{\bs{i},n}^{\pm 1}$,
 are contained in $U_q^{\A}(\g)\tensor U_{q,e}^{\A[2Q_0^-]}(\gt)$.
A similar argument works for $U_{q,\chi}^k(\tg)$.
\end{proof}

Next we prove a PBW-type result for $U_{q,\chi}^k(\gt)$.
Let us define $\tE_{\bs{i},\alpha}, \tK_{\mu}, \tF_{\bs{i},\alpha} \in U_q^{\A[P]}(\g)$ as
$(q_{\alpha} - q_{\alpha}^{-1})\aE_{\bs{i},\alpha}, e_{-\mu}K_{\mu}, (q_{\alpha} - q_{\alpha}^{-1})\aF_{\bs{i},\alpha}$ respectively.

\begin{lemm} \label{lemm:right system case}
Set $\cl{U} := \mathrm{span}_{\A[2Q^-]}\{\tF_{\bs{i}}^{\Lambda^-}\tK_{\mu}\tE_{\bs{i}}^{\Lambda^+}\}_{(\Lambda^{\pm},\mu)} \subset U_q^{\A[P]}(\g)$.
Then this is closed under multiplication.
Moreover, for all $1 \le l \le N$, 
it coincides with the $\A[2Q^-]$-linear span of elements of
the form 
$\tE_{\bs{i}}^{\Lambda_2^+}\tF_{\bs{i}}^{\Lambda^-}\tK_{\mu}\tE_{\bs{i}}^{\Lambda_1^+}$ where $\Lambda_1^+ \le l < \Lambda_2^+$.
\end{lemm}
\begin{proof}
It is not difficult to see that
$\cl{U}^- := \mathrm{span}_{\A[2Q^-]}\{\tF_{\bs{i}}^{\Lambda^-}\tK_{\mu}\}$ and
$\cl{U}^+ := \mathrm{span}_{\A[2Q^-]}\{\tK_{\mu}\tE_{\bs{i}}^{\Lambda^+}\}$
are closed under multiplication. To see
$xy \in \cl{U}$ for $x \in \cl{U}^+$ and $y \in \cl{U}^-$,
one should note that the argument in the proof of Proposition \ref{prop:mixed LS relation} (iii) works for computing the expansion of $\tE_{\bs{i}}^{\Lambda^+}\tF_{\bs{i}}^{\Lambda^-}$. Hence this product can be expressed as an
$\A$-linear combination of elements of the form $\tF_{\bs{i}}^{\Gamma^-}K_{2\mu}^{-1}\tE_{\bs{i}}^{\Gamma^+}$ with $\mu \in Q^+$,
which shows the required property since $K_{2\mu}^{-1} = e_{-2\mu}\tK_{2\mu}^{-1}$.

Next we show the latter half of the statement
by downward induction on $l$. If $l = N$,
there is nothing to prove. For the induction step,
assume the statement holds for a fixed $l$. It suffices to show that
$\cl{U}' := \mathrm{span}_{\A[2Q^-]}\{\tE_{\bs{i}}^{\Lambda_2^+}\tF_{\bs{i}}^{\Lambda^-}\tK_{\mu}\tE_{\bs{i}}^{\Lambda_1^+}\mid\Lambda_1^+ < l \le \Lambda_2^+\}$ is closed under multiplication,
which, combining with the induction hypothesis, implies that the linear span is the
$\A[2Q^-]$-subalgebra generated by
$(\tF_{\bs{i},\alpha})_{\alpha \in \roots^+}, (\tK_{\mu})_{\mu \in P}$ and $(\tE_{\bs{i},\alpha})_{\alpha \in\roots^+}$,
which is nothing but $\cl{U}$.

We only prove that $\cl{U}'$ is invariant under
the right multiplication by $\tE_{\bs{i},l}$
since this fact and the induction hypothesis imply
the invariance for other elements.

Set $\cl{U}_l^{\pm}$ as follows:
\begin{align*}
 \cl{U}_l^- &= \aspan_{\A[2Q^-]}\{\tF_{\bs{i}}^{\Lambda^-}\tK_{\mu}\tE^{\Lambda^+}\mid \mu \in P, \Lambda^+ < l < \Lambda^-\} \\
 \cl{U}_l^+ &= \aspan_{\A[2Q^-]}\{\tE_{\bs{i}}^{\Lambda^+}\tK_{\mu}\tF_{\bs{i}}^{\Lambda^-}\mid \mu \in P, \Lambda^- < l < \Lambda^+\}.
\end{align*}
These are closed under multiplication by Proposition \ref{prop:mixed LS relation} (i).
Also note that $\cl{U}'$ is spanned by elements of the form $\tE_{\bs{i},l}^px^+\tF_{\bs{i},l}^qx^-$ with $x^{\pm} \in \cl{U}_l^{\pm}$ and
we have
\begin{align*}
 \tF_{\bs{i},l}\tE_{\bs{i},l}
= q_{\alpha^{\bs{i}}_l}^{-2}\tE_{\bs{i},l}\tF_{\bs{i},l} - (q_{\alpha^{\bs{i}}_l} - q_{\alpha^{\bs{i}}_l}^{-1})(1 - e_{-2\alpha^{\bs{i}}_l}\tK_{-2\alpha^{\bs{i}}_l}).
\end{align*}
Therefore it suffices to show that $\cl{U}_l^{\pm}\tE_{\bs{i},l} \subset \tE_{\bs{i},l}\cl{U}_l^{\pm} + \cl{U}_l^{\pm}$. At first take $x^- \in \cl{U}_l^-$.
By using the triangular decomposition and Proposition \ref{prop:LS relation},
we may assume $x^- = \tF_{\bs{i}}^{\Lambda^-}$ with $\Lambda^- > l$.
If $\Lambda^- = \delta_m$ for some $m > l$,
the mixed Levend\"{o}rskii-Soibelman relation (\ref{eq:mixed LS relation})
proves $\tF_{\bs{i},m}\tE_{\bs{i},l} - q^{-(\alpha^{\bs{i}}_l,\alpha^{\bs{i}}_m)}\tE_{\bs{i},l}\tF_{\bs{i},m} \in \cl{U}_l^-$.
Then we have
\begin{align*}
 \tF_{\bs{i},m_1}\tF_{\bs{i},m_2}\cdots \tF_{\bs{i},m_n}\tE_{\bs{i},l}
- q^{-(\alpha^{\bs{i}}_l, \alpha^{\bs{i}}_{m_1} + \alpha^{\bs{i}}_{m_2} + \cdots + \alpha^{\bs{i}}_{m_l})}\tE_{\bs{i},l}\tF_{\bs{i},m_1}\tF_{\bs{i},m_2}\cdots \tF_{\bs{i},m_n}
\in \cl{U}_l^-
\end{align*}
for any sequence $(m_1,m_2,\cdots,m_n)$ with $m_j > l$, which is
nothing but the required property.
A similar argument works for $x^+ \in \cl{U}_l^+$.
\end{proof}

\begin{prop} \label{prop:PBW theorem for gt}
Fix a reduced expression $s_{\bs{i}}$ of the longest element compatible with $\roots_0^+$.
Then $\{\hF_{\bs{i}}^{(\Lambda^-)}\hK_{\mu}\hE_{\bs{i}}^{\Lambda^+}\}_{(\Lambda^{\pm},\mu)}$ is a basis of $U_{q,\chi}^k(\gt)$.
\end{prop}
\begin{proof}
We may assume $k = \A[2Q_0^-]$.
 Set $l := \abs{\roots^+\cap \roots_0^+}$.
Then there is $w \in W$ and a reduced expression $s_{\bs{j}}$ of $w_0$ as Lemma \ref{lemm:rotation lemma} for $s_{\bs{i}}$ and 
$N - l$.
Let $t_w$ be an automorphism of $\A[P]$ which sends $t_w(e_{\lambda}) = e_{w(\lambda)}$ and $\phi$ be an automorphism on $U_q^{\A[P]}$ defined as $t_w\tensor \cl{T}_w$. Then
$\phi(\cl{U})$ is closed under multiplication and spanned by elements
of the form
$\cl{T}_w(\tE_{\bs{j}})^{\Lambda_2^+}\cl{T}_w(\tF_{\bs{j}})^{\Lambda^-}\tK_{\mu}\cl{T}_w(\tE_{\bs{j}})^{\Lambda_1^+}$ with $\Lambda_1^+ \le l < \Lambda_2^+$. Hence Lemma \ref{lemm:right system case} implies
that $\mathrm{span}_k \{\vF_{\bs{i}}^{\Lambda^-}\hK_{\mu}\hE_{\bs{i}}^{\Lambda^+}\}_{(\Lambda^{\pm},\mu)}$ is closed under multiplication.

Now fix $(\Lambda_1^{\pm},\mu_1)$ and $(\Lambda_2^{\pm},\mu_2)$ and
consider the following two expansions with 
$C_{\Lambda^{\pm},\mu} \in \A$:
\begin{align*}
 (\aF_{\bs{i}}^{(\Lambda_1^-)}K_{\mu_1}\tE_{\bs{i}}^{\Lambda_1^+})(\aF_{\bs{i}}^{(\Lambda_2^-)}K_{\mu_2}\tE_{\bs{i}}^{\Lambda_2^+})
&= \sum_{\Lambda^{\pm},\mu} C_{\Lambda^{\pm},\mu}\aF_{\bs{i}}^{(\Lambda^-)}K_{\mu}\tE^{\Lambda^+}.
\end{align*}
Then there is $\nu(\Lambda^{\pm},\mu) \in P$ for each
$(\Lambda^{\pm},\mu)$ with which we have
\begin{align*}
 (\hF_{\bs{i}}^{(\Lambda_1^-)}\hK_{\mu_1}\hE_{\bs{i}}^{\Lambda_1^+})(\hF_{\bs{i}}^{(\Lambda_2^-)}\hK_{\mu_2}\hE_{\bs{i}}^{\Lambda_2^+})
= \sum_{\Lambda^{\pm},\mu} e_{\nu(\Lambda^{\pm},\mu)}C_{\Lambda^{\pm},\mu}\hF_{\bs{i}}^{(\Lambda^-)}\hK_{\mu}\hE_{\bs{i}}^{\Lambda^+}.
\end{align*}
On the other hand, this weight appears also in the following
comparison:
\begin{align*}
 (\tF_{\bs{i}}^{\Lambda_1^-}K_{\mu_1}\tE_{\bs{i}}^{\Lambda_1^+})(\tF_{\bs{i}}^{\Lambda_2^-}K_{\mu_2}\tE_{\bs{i}}^{\Lambda_2^+})
&= \sum_{\Lambda^{\pm},\mu} c_{\Lambda^{\pm},\mu}\tF_{\bs{i}}^{\Lambda^-}K_{\mu}\tE^{\Lambda^+}, \\
 (\vF_{\bs{i}}^{\Lambda_1^-}\hK_{\mu_1}\hE_{\bs{i}}^{\Lambda_1^+})(\vF_{\bs{i}}^{\Lambda_2^-}\hK_{\mu_2}\hE_{\bs{i}}^{\Lambda_2^+})
&= \sum_{\Lambda^{\pm},\mu} e_{(\Lambda^{\pm},\mu)}c_{\Lambda^{\pm},\mu}\vF_{\bs{i}}^{\Lambda^-}\hK_{\mu}\hE_{\bs{i}}^{\Lambda^+}.
\end{align*}
Then the discussion above shows that
each $\nu(\Lambda^{\pm},\mu)$ is contained
in the image of $t_w(2Q^-) = 2Q_0^-$, which completes the proof.
\end{proof}

A similar argument works to prove the following.

\begin{prop} \label{prop:PBW for tg}
Fix a reduced expression $s_{\bs{i}}$ of the longest element compatible with $\roots_0^+$.
Then $\{\vF_{\bs{i}}^{\Lambda^-}\vK_{\mu}\vE_{\bs{i}}^{(\Lambda^+)}\}_{(\Lambda^{\pm},\mu)}$ is a basis of $U_q^k(\tg)$.
\end{prop}

\subsection{Twisted parabolic induction via deformed QEA}
\label{subsection:twisted parabolic induction via deformed QEA}
At the end of this section, we define twisted parabolic induction,
which play an important role to construct quantum semisimple coadjoint orbits.

Let $S$ be a set of simple roots and $\roots_0^+$
be a positive system of $\roots$ containing $\roots_S^+$.
Fix a reduced expression $s_{\bs{i}}$ compatible with $\roots_0^+$.
We also assume that $\roots_S^+ = \{\alpha^{\bs{i}}_k\}_{k = N - N_S + 1}^N$ holds:
\begin{align*}
 \underbrace{\alpha^{\bs{i}}_1,\, \alpha^{\bs{i}}_2,\,\dots,\, \alpha^{\bs{i}}_{l}}_{\roots^+\setminus\roots_0^+}, \,\underbrace{\alpha^{\bs{i}}_{l + 1},\, \alpha^{\bs{i}}_{l + 2},\, \dots, \alpha^{\bs{i}}_{N - N_S}}_{\roots_0^+\setminus\roots_S^+},\, \underbrace{\alpha^{\bs{i}}_{N - N_S + 1},\, \alpha^{\bs{i}}_{N - N_S + 2}\dots,\, \alpha^{\bs{i}}_N}_{\roots_S^+}.
\end{align*}
where $l = \abs{\roots^+\setminus\roots_0^+}$.
For an $\A[2Q_0^-]$-algebra $(k,\chi)$ with
$\chi_{-2\alpha} = 1$ for $\alpha \in \roots_S^+$,
we introduce the following $k$-subalgebras:
\begin{itemize}
 \item $U_{q,\chi}^{k}(\frpt)$, generated by
$(\hK_{\mu})_{\mu \in P},\,(\hE_{\bs{i},\alpha})_{\alpha\in \roots^+}$ and $(\hF_{\bs{i},\alpha}^{(r)})_{\alpha \in \roots_S^+, r\ge 0}$.
 \item $U_{q,\chi}^k(\lt)$, generated by
$(\hK_{\mu})_{\mu \in P},\,(\hE_{\bs{i},\alpha})_{\alpha\in \roots_S^+}$ and $(\hF_{\bs{i},\alpha}^{(r)})_{\alpha \in \roots_S^+,r \ge 0}$.
 \item $U_{q,\chi}^k(\utp)$, generated by $(\hE_{\bs{i},\alpha})_{\alpha \in \roots^+\setminus\roots_S^+}$.
 \item $U_{q,\chi}^k(\utm)$, generated by $(\hF_{\bs{i},\alpha}^{(r)})_{\alpha \in \roots^+\setminus\roots_S^+, r\ge 0}$.
\end{itemize}

Due to the additional requirement on $s_{\bs{i}}$, these do not
depend on the choice of such a reduced expression. Moreover these
are spanned by the PBW basis introduced in Proposition \ref{prop:PBW theorem for gt}.
Hence the multiplication map induces the following tensor product
decompositions:
\begin{align}
U_{q,\chi}^k(\frpt)&\cong U_{q,\chi}^k(\lt)\tensor U_{q,\chi}^k(\utp),\\ \label{eq:tensor product decomposition of gt for verma}
U_{q,\chi}^k(\gt)&\cong U_{q,\chi}^k(\utm)\tensor U_{q,\chi}^k(\frpt).
\end{align}

For $\chi \equiv 1$, these algebras are simply denoted by
$U_q^k(\gt), U_q^k(\frpt), U_q^k(\lt)$ and $U_q^k(\fru_S^{\pm}{}^{\sim})$ respectively.
These can be thought as subalgebras of $U_q^k(\g)$ by 
the identifications $\hF_{\bs{i},\alpha}^{(r)} = F_{\bs{i},\alpha}^{(r)}, \hK_{\mu} = K_{\mu}, \hE_{\bs{i},\alpha} = \tE_{\bs{i},\alpha}$.

\begin{lemm}
There is a canonical isomorphism from $U_{q,\chi}^k(\frpt)$
to $U_q^k(\frpt)$,
which is defined as follows on the generators:
\begin{align*}
 \hF_{\bs{i},\alpha}^{(r)}\longmapsto F_{\bs{i},\alpha}^{(r)},\quad \hK_{\mu}\longmapsto K_{\mu},\quad \hE_{\bs{i},\alpha}\longmapsto \tE_{\bs{i},\alpha}.
\end{align*}
\end{lemm}
\begin{rema}
In light of this fact,
we identify $U_{q,\chi}^k(\frpt)$ with $U_q^k(\frpt)$.
\end{rema}
\begin{proof}
Using Proposition \ref{prop:PBW theorem for gt} we can define
a $k$-linear map from $U_{q,\chi}^k(\frpt)$ to $U_q^k(\frpt)$ sending
$\hF_{\bs{i}}^{(\Lambda^-)}\hK_{\mu}\hE_{\bs{i}}^{\Lambda^+}$
to $F_{\bs{i}}^{(\Lambda^-)}K_{\mu}\tE_{\bs{i}}^{\Lambda^+}$.
What we have to show is multiplicativity of this map.
By Proposition \ref{prop:mixed LS relation},
$\tE_{\bs{i}}^{\Lambda_0^+}F_{\bs{i}}^{(\Lambda_0^-)}$ with $\Lambda^+ \le N - N_S < \Lambda^-$ is a linear combination of $F_{\bs{i}}^{(\Lambda^-)}\tE_{\bs{i}}^{\Lambda^+}$ with $\Lambda_0^+ \le N - N_S < \Lambda_0^-$. Hence we have
\begin{align*}
 F_{\bs{i}}^{(\Lambda_1^-)}K_{\mu_1}\tE_{\bs{i}}^{\Lambda_1^+}
F_{\bs{i}}^{(\Lambda_1^-)}K_{\mu_1}\tE_{\bs{i}}^{\Lambda_1^+}
= \sum_{\substack{\Lambda^{\pm},\mu \in Q_S\\ \Lambda^- > N - N_S}}
C_{\Lambda^{\pm},\mu}F_{\bs{i}}^{(\Lambda^-)}K_{\mu_1 + \mu_2 + \mu}\tE_{\bs{i}}^{\Lambda^+}
\end{align*}
when $\Lambda_i^- > N - N_S$. Now the multiplicativity follows from
this expansion and the definition of the generators
$\hF_{\bs{i},\alpha}^{(r)}, \hK_{\mu}, \hE_{\bs{i},\alpha}$.
\end{proof}
Note that there is a canonical homomorphism $\map{\pi_S}{U_q^k(\frpt)}{U_q^k(\lt)}$ which is identical on $U_q^k(\lt)$ and equal to $\eps$
on $U_q^k(\utp)$. All $U_q^k(\lt)$-modules are regarded as
$U_q^k(\frpt)$-modules via this homomorphism.
\begin{defn}
We define a functor
$\map{\dqind}{\mod{U_q^k(\lt)}}{\mod{U_{q,\chi}^k(\gt)}}$
as $U_{q,\chi}^k(\gt)\tensor_{U_q^k(\frpt)}\tend$.
For a $U_{q,\chi}^k(\gt)$-module $V$, $\dqind{V}$
is called the $(\roots_0^+,\chi)$-twisted
parabolic induction module of $V$.
\end{defn}
\begin{rema}
By the tensor product decomposition (\ref{eq:tensor product decomposition of gt for verma}),
$\dqind{V}$ can be
naturally identified with $U_{q,\chi}^k(\utm)\tensor V$
as $U_{q,\chi}^k(\utm)$-modules.
\end{rema}

Let $M$ be a $U_{q,\chi}^k(\gt)$-module. We say that $m \in M$ is
\emph{$S$-maximal} when $xm = \pi_S(x)m$ for $x \in U_q^k(\frpt)$,
which is the same thing as saying that $\hE_{\bs{i},\beta}m = 0$
for $\beta \in \roots^+\setminus\roots_S^+$.
The following universal property of $(\roots_0^+,\chi)$-twisted
parabolic induction is proven in an analogous way to the usual setting.
\begin{lemm} \label{lemm:universal property}
Let $M$ be a $U_q^k(\gt)$-module. Then $M_{S\text{-}\max}$, the set of
$S$-maximal vectors in $M$, forms a $U_q^k(\lt)$-module.
Moreover we have a canonical isomorphism:
\begin{align*}
\Hom_{U_q^k(\lt)}(V,M_{S\text{-}\max})
\cong \Hom_{U_{q,\chi}^k(\gt)}(\dqind{V}, M).
\end{align*}
\end{lemm}

Let $\map{\qind}{\mod{U_q^k(\lt)}}{\mod{U_q^k(\gt)}}$ be
the usual parabolic induction $U_q^k(\gt)\tensor_{U_q^k(\frpt)}\tend$.
In order to compare the $(\roots_0^+,\chi)$-twisted parabolic
induction with the usual parabolic induction, we consider subalgebras
$U_{q,\chi}^k(\gt)_{2Q}$ defined by restricting the Cartan part $\aspan_k\{\hK_{\mu}\}_{\mu \in P}$ to $\aspan_k\{\hK_{2\mu}\}_{\mu \in Q}$
and similar variants for the other subalgebras. 
Note that $\dqind V$ is
naturally isomorphic to $U_{q,\chi}^k(\gt)_{2Q}\tensor_{U_{q,\chi}^k(\frpt)_{2Q}} V$ as $U_{q,\chi}^k(\gt)_{2Q}$-module and the same thing
can be said for $\qind{V}$.

\begin{lemm} \label{lemm:comparison with usual verma}
Assume that $\map{\chi}{2Q^-}{k}$ extends to a group homomorphism $\map{\chi}{2Q}{k^{\times}}$ with $\chi|_{2Q_S} \equiv 1$.
Let $k_{\chi}$ be a $U_q^k(\lt)$-module
which is $k$ as a $k$-module and
$\hF_{\bs{i},\alpha}^{(r)}, \hK_{2\mu}, \hE_{\bs{i},\beta}$ act on
it as $0,\chi_{2\mu}, 0$ respectively.
\begin{enumerate}
 \item There is an isomorphism $U_{q,\chi}^k(\gt)_{2Q}\cong U_q^k(\gt)_{2Q}$ with:
\begin{align*}
\hF_{\bs{i},\alpha}^{(r)}\longmapsto 
\begin{cases}
 \chi_{2\alpha}^rF_{\bs{i},\alpha}^{(r)} &(\alpha \in \roots^+\setminus\roots_0^+) \\
 F_{\bs{i},\alpha}^{(r)}& (\alpha \in \roots^+\cap\roots_0^+)
\end{cases},\quad
\hK_{2\mu}\longmapsto \chi_{-2\mu}K_{2\mu},\quad
\hE_{\bs{i},\alpha}\longmapsto \tE_{\bs{i},\alpha}.
\end{align*}
 \item Let $V$ be a $U_q^k(\lt)$-module. Then $\dqind{V}$ is
canonically isomorphic to $\qind{(V\tensor k_{\chi})}$
as a $U_q^k(\gt)_{2Q}$-module.
\end{enumerate}
\end{lemm}
\begin{proof}
(i) This can be immediately verified by the definition of $U_{q,\chi}^k(\gt)$.

(ii) This can be seen from the universal property of the $(\roots_0^+,\chi)$-twisted parabolic induction and that of the usual parabolic induction.
\end{proof}

\section{$S$-maximal vectors in the formal setting}
\label{section:S-maximal vectors in the formal setting}
\subsection{Notation}
In Section \ref{section:S-maximal vectors in the formal setting} and
\ref{section:globalization of the 2-cocycle},
we work on $(k\bbh, s = e^{h/2L})$
where $k$ is a $\Q$-algebra. In this case it
is convenient to consider $h$-adically completed 
tensor products, which is denoted by $\tend\htensor \tend$ in this paper. 

In the formal setting we consider the $h$-adic Drinfeld-Jimbo
deformation $U_h^k(\g) := k\bbh\htensor U_h(\g)$,
where the definition of $U_h(\g)$ is described in
\cite[Section 6.1.3, Definition 2]{MR1492989} for instance.
This is isomorphic to $(k\tensor U(\g))\bbh$ as topological
$k\bbh$-algebra. Moreover we can see that a topological
$U_h^k(\g)$-module admitting a finite $k\bbh$-basis
is isomorphic to $(k\tensor V)\bbh$ where $V$ is
a finite dimensional representation of $\g$ (\cite[Section 7.1.3, Proposition 10]{MR1492989}).
Such a module is called a \emph{finite integrable} module
in this paper and the category of finite integrable modules is denoted by
$\rep_h^{\fin}\g$.

Fix $S, \roots_0^+, s_{\bs{i}}$ and $\chi$ as those in Subsection \ref{subsection:twisted parabolic induction via deformed QEA}.
In the formal setting it is rather appropriate to add
$\hH_{\alpha}$ to $U_{q,\chi}^{k\bbh}(\gt)$
so that $\hK_{\alpha} = \exp{hd_{\alpha}\hH_{\alpha}}$.
The precise definition is the following.

\begin{defn}
Set $U^k(\h) := k\tensor U(\h)$ and $U_h^k(\h)$ be the $h$-adic
completion of $U^k(\h)\tensor k\bbh$. This is equipped with
a canonical coproduct such that $\Delta(H) = H\tensor 1 + 1\tensor H$
for $H \in \h$.
Consider
$U_{q,\chi}^{k\bbh}(\gt)\tensor_{k\bbh} U_h^k(\h)$ with
the following multiplication:
\[
 (x \tensor f)(y\tensor g) = (\wt y)(f_{(1)})(xy \tensor f_{(2)}g),
\]
where $\wt y$, which is an element of $\h^*$, is considered as a character on $U_h^k(\h)$.
Then We define $U_{h,\chi}^k(\gt)$ as a quotient of
$U_{q,\chi}^{k\bbh}(\gt)\tensor_{k\bbh} U_h^k(\h)$ by the ideal
generated by $\hK_{\lambda}\tensor 1 - 1\tensor \exp(h\lambda)$,
where $\lambda$ is considered as an element of $\h$ uniquely determined by $\alpha(\lambda) = (\alpha,\lambda)$ for $\alpha \in \h^*$.
\end{defn}

We refer definitions and properties developed in
Section \ref{sect:deformed QEA} to apply them to $U_{h,\chi}^k(\gt)$
in a suitably modified form. For instance $U_{h,\chi}^k(\gt)$
has a canonical left coaction of $U_h^k(\g)$, the $h$-adic
Drinfeld-Jimbo deformation of $\g$.

\subsection{Generalized maximal vector}
Though what we would like to show in this section is a kind of
semisimplicity of the tensor product of
finite integrable $U_h^k(\g)$-modules and
$(\roots_0^+,\chi)$-twisted parabolically induced modules,
we show a stronger result: there is
an operator in $U_h^k(\g\times \frl_S)$
which produces $S$-maximal vectors in the tensor product modules.
To discuss such an operator, we introduce
the following generalization of $S$-maximal vectors.
Like the Sweedler notation, $x_{(\fru_S^-)}\tensor x_{(\frp_S)}$
denotes the image of $x \in U_{h,\chi}^k(\gt)$ under the isomorphism
$U_{h,\chi}^k(\gt)\cong U_{h,\chi}^k(\utm)\tensor U_{h,\chi}^k(\frpt)$. Similarly $x \in W\htensor V$ is denoted by $x_{(W)}\tensor x_{(V)}$. For instance, the image of $x$ under a linear mapping $w\tensor v \longmapsto w\tensor a \tensor v$ is denoted by $x_{(W)}\tensor a \tensor x_{(V)}$.

\begin{defn} \label{defn:generalized maximal vector}
Let $V$ be a complete $U_h^k(\frl_S)$-module and $W$ be a complete
$U_h^k(\g)$-module.
We say that $(m_{\Lambda})_{\Lambda \le N - N_S} \subset W\htensor V$ is
a \emph{generalized $S$-maximal vector} if it satisfies
the following equation for all $\beta \in\roots^+\setminus\roots_S^+$:
\begin{align} \label{eq:recursion formula}
 \sum_{\Lambda} \hE_{\bs{i},\beta,(1)}m_{\Lambda,(W)}\tensor (\hE_{\bs{i},\beta,(2)}\hF_{\bs{i}}^{(\Lambda)})_{(\fru_S^-)}\tensor \pi_S((\hE_{\bs{i},\beta,(2)}\hF_{\bs{i}}^{(\Lambda)})_{(\frp_S)})m_{\Lambda,(V)}
= 0. 
\end{align}
\end{defn}
\begin{rema}
Strictly speaking, these equalities should be interpreted after
expanding $(\hE_{\bs{i},\beta,(2)}\hF_{\bs{i}}^{(\Lambda)})_{(\fru_S^-)}$ along the PBW basis. Then we obtain well-defined
equations since, for any fixed $\Gamma$,
only finitely many terms produce $\hF_{\bs{i}}^{(\Gamma)}$.
\end{rema}

In the rest of this subsection, we assume
$\chi_{-2\alpha} = 1$ for $\alpha \in \roots_S^+$ and $\chi_{-2\alpha} - 1 \in k\bbh^{\times}$
for $\alpha \in \roots_0^+\setminus\roots_S^+$.

\begin{prop} \label{prop:uniqueness of generalized maximal vector}
Let $V$ (resp. W) be a complete $U_h^k(\frl_S)$-module
(resp. $U_h^k(\g)$-module). Then, for any $m \in W\htensor V$,
a generalized $S$-maximal vector $(m_{\Lambda})_{\Lambda}$ with
$m_0 = m$ exists and is unique.
\end{prop}

To prove this proposition, we consider a generalized version.
Let $(p_{\Gamma})_{\Gamma \le N - N_S}$ be the dual basis
of the PBW basis $(\hF_{\bs{i}}^{(\Gamma)})_{\Gamma \le N - N_S}$
of $U_{h,\chi}^k(\utm)$. Then, for $x \in W\htensor U_{h,\chi}^k(\utm)\htensor V$, we define the coefficient of $\hF_{\bs{i}}^{(\Gamma)}$
in $x$ as $(\id\tensor p_{\Gamma}\tensor \id)(x) \in W\htensor V$.

We say that $(m_{\Lambda})_{\Lambda \le N - N_S}$ is a
\emph{weak solution} of (\ref{eq:recursion formula}) for $\beta = \alpha^{\bs{i}}_l$ if the coefficient of $\hF_{\bs{i}}^{(\Gamma)}$ in the RHS
is $0$ for all $\Gamma$ with $\Gamma \ge l.$

\begin{lemm} \label{lemm:uniqueness and existence of partly maximal vector}
Fix a positive integer $ 0 \le l \le N - N_S$.
\begin{enumerate}
 \item If $(m_{\Lambda})_{l < \Lambda \le N - N_S}$ is given,
this extends a unique family $(m_{\Lambda})_{\Lambda \le N - N_S}$
which is a weak solution of (\ref{eq:recursion formula}) for all $\alpha^{\bs{i}}_m$ with $1 \le m \le l$.
 \item If $(m_{\Lambda})_{\Lambda \le N - N_S}$ is a weak solution
of (\ref{eq:recursion formula}) for all $\alpha^{\bs{i}}_{m}$
with $1 \le m \le l$, it actually satisfies (\ref{eq:recursion formula}) for $\alpha^{\bs{i}}_m$ with $1 \le m \le l$.
\end{enumerate}
\end{lemm}
Let $\map{\Pi}{Q}{Q_{\simples\setminus S}}$ be the projection with respect to a basis $\Delta$.
\begin{proof}
(i) Fix a multi-index $\Gamma$ with $l \not<\Gamma \le \abs{\roots^+\setminus\roots_S^+}$
and consider a decomposition $\Gamma = \Gamma' + \delta_i$
such that $\gamma_1 = \gamma_2 = \cdots = \gamma_{i - 1} = 0$ and $\gamma_i = \gamma'_i + 1$. In this case $i \le l$ holds.

We look at the coefficient of
$\hF_{\bs{i}}^{(\Gamma')}$ in the RHS of (\ref{eq:recursion formula})
with $\beta = \alpha^{\bs{i}}_i$.
Note that $\Delta(\hE_{\bs{i},i})$ is
a linear combination of $K_{\Lambda^r\cdot\alpha^{\bs{i}}}^{-1}\tE_{\bs{i}}^{\Lambda^l}\tensor \hE_{\bs{i}}^{\Lambda^r}$ with $\Lambda^r \le i \le \Lambda^l$ and $\Lambda^l\cdot\alpha^{\bs{i}} + \Lambda^r\cdot\alpha^{\bs{i}} = \alpha^{\bs{i}}_i$ as we can see in
(\ref{eq:coproduct of E}).

If $\Pi(\Lambda^l\cdot\alpha^{\bs{i}})\neq 0$ holds, $p_{\Gamma'}((\hE_{\bs{i}}^{\Lambda^r}\hF_{\bs{i}}^{(\Lambda)})_{(\fru_S^-)})\pi_S((\hE_{\bs{i}}^{\Lambda^r}\hF_{\bs{i}}^{(\Lambda)})_{(\frp_S)}) \neq 0$ implies 
$\Pi(\Lambda\cdot\alpha^{\bs{i}}) = \Pi(\Gamma'\cdot\alpha^{\bs{i}}) + \Pi(\Lambda^r\cdot\alpha^{\bs{i}}) < \Pi(\Gamma\cdot\alpha^{\bs{i}})$.

If $\Pi(\Lambda^l\cdot\alpha^{\bs{i}}) =0$ and $\lambda_1 = \lambda_2 = \cdots = \lambda_{i - 1} = 0$ holds, $\Pi(\Lambda\cdot\alpha^{\bs{i}}) = \Pi(\Gamma\cdot\alpha^{\bs{i}})$ is necessary for
$p_{\Gamma'}((\hE_{\bs{i}}^{\Lambda^r}\hF_{\bs{i}}^{(\Lambda)})_{(\fru_S^-)})\pi_S((\hE_{\bs{i}}^{\Lambda^r}\hF_{\bs{i}}^{(\Lambda)})_{(\frp_S)}) \neq 0$ by the discussion above.
Moreover $\Lambda^r < i$ holds when $\Lambda^l \neq 0$, since
$\Lambda^l\cdot\alpha^{\bs{i}} + \Lambda^r\cdot\alpha^{\bs{i}} = \alpha^{\bs{i}}_i$.
Then Proposition \ref{prop:mixed LS relation} (ii) implies $p_{\Gamma'}((\hE_{\bs{i}}^{\Lambda^r}\hF_{\bs{i}}^{(\Lambda)})_{(\fru_S^-)})\pi_S((\hE_{\bs{i}}^{\Lambda^r}\hF_{\bs{i}}^{(\Lambda)})_{(\frp_S)}) \in h^{\abs{\Lambda^l}}U_{h}^k(\lt)$.
On the other hand $\Lambda^r = \delta_i$ holds when $\Lambda^l = 0$,
since $i \le \Lambda^r$ and $\Lambda'\cdot\alpha^{\bs{i}} \neq \alpha^{\bs{i}}_i$ when $i < \Lambda'$ by \cite[Theorem, p.662]{MR1169886}.
Hence we have to look at the coefficient of $\hF_{\bs{i}}^{(\Gamma')}$in the following summation:
\begin{align*}
\sum_{\Lambda} K_{\beta}^{-1}m_{\Lambda,(W)}\tensor [\hE_{\bs{i},\beta}, \hF_{\bs{i}}^{(\Lambda)}]_{q,(\fru_S^-)}\tensor \pi_S([\hE_{\bs{i},\beta}, \hF_{\bs{i}}^{(\Lambda)}]_{q,(\frp_S)})m_{\Lambda,(V)}.
\end{align*}
Consider a decomposition $\Lambda = \lambda_i\delta_i + \Lambda_{>i}$ with $\Lambda_{>i} > i$. Then we have
\begin{align*}
 [\hE_{\bs{i},\beta}, \hF_{\bs{i}}^{(\Lambda)}]_q
= [\hE_{\bs{i},\beta},\hF_{\bs{i},\beta}^{(\lambda_i)}]_q\hF_{\bs{i}}^{(\Lambda_{>i})}
 + q_{\beta}^{2\lambda_i}\hF_{\bs{i},i}^{(\lambda_i)}[\hE_{\bs{i},\beta},\hF_{\bs{i}}^{(\Lambda_{>i})}]_q.
\end{align*}
When $\Lambda = \Gamma$,
the first term produces $\hF_{\bs{i}}^{(\Gamma')}$ and
the second term does not.
Moreover we can determine its coefficient completely since we have
\begin{align*}
[\hE_{\bs{i},\beta},\hF_{\bs{i},\beta}^{(r)}]_q =
\begin{cases}
 q_{\beta}^{2r}\hF_{\bs{i},\beta}^{(r - 1)}(q_{\beta}^{1 - r} - q_{\beta}^{r - 1}e_{-2\beta}\hK_{-2\beta}) & (\beta \in \roots_0^+), \\
 q_{\beta}^{2r}\hF_{\bs{i},\beta}^{(r - 1)}(q_{\beta}^{1 - r}e_{2\beta} - q_{\beta}^{r - 1}\hK_{-2\beta}) & (\beta \not\in \roots_0^+).
\end{cases}
\end{align*}
Also note that the coefficient of $K_{\beta}^{-1}\tensor \hE_{\bs{i},\beta}$ in $\Delta(\hE_{\bs{i},\beta})$ is $1$ by (\ref{eq:coproduct of E}).

When $\Lambda \neq \Gamma$,
the coefficient of $\hF_{\bs{i}}^{(\Gamma')}$ in $[\hE_{\bs{i},\beta},\hF_{\bs{i},\beta}^{(\lambda_i)}]_q\hF_{\bs{i}}^{(\Lambda_{>i})}$
is $0$ and that in $\hF_{\bs{i},i}^{(\lambda_i)}[\hE_{\bs{i},\beta},\hF_{\bs{i}}^{(\Lambda_{>i})}]_q$ is in $hU_{h}^k(\frl_S)$
by Proposition \ref{prop:mixed LS relation} (i).

Take $\nu \in Q_{\simples\setminus S}^+$ and
set $B_{l,\nu} := \{\Gamma\mid \Pi(\Gamma\cdot\alpha^{\bs{i}}) = \nu,\,l \not<\Gamma \le \abs{\roots^+\setminus\roots_S^+}\}$. Then $n_{l,\nu} := \abs{B_{l,\nu}} < \infty$. We also consider an enumeration $\{\Gamma_i\}_i$ on $B_{l,\nu}$ along the lexicographic order.
Then the consideration above implies that we can find
a $n_{l,\nu}\times n_{l,\nu}$-matrix $A_{l,\nu}$
and $(u_i)_i \in (W\htensor V)^{n_{l,\nu}}$ such that
\begin{itemize}
 \item the $i$-th entry of $A_{l,\nu}(m_{\Gamma_i})_i + (u_i)_i$
is the coefficient of $\hF_{\bs{i}}^{(\Gamma'_i)}$ in the LHS of (\ref{eq:recursion formula}),
 \item each entry of $A_{l,\nu}$ is in $U_h^k(\g\times \frl_S)$,
 \item each $u_i$ is in $U_h^k(\g\times \frl_S)$-submodule generated
by $m_{\Lambda}$ with $\Pi(\Lambda\cdot\alpha^{\bs{i}}) < \nu$,
 \item the matrix $A_{l,\nu}$ is of the following form modulo $h$:
\begin{align*}
\begin{pmatrix}
 \sigma_1(1 - e_{-2\sigma_1\beta_1}) & \ast & \cdots & \ast \\
 0 & \sigma_2(1 - e_{-2\sigma_2\beta_2}) & \cdots & \ast \\
\vdots & \vdots & \ddots & \vdots \\
 0 & 0 & \cdots & \sigma_{n_{l,\nu}}(1 - e_{-2\sigma_{n_{l,\nu}}\beta_{n_{l,\nu}}})
\end{pmatrix},
\end{align*}
where $\sigma_i = 1$ when $\beta_i \in \roots^+\cap \roots_0^+$
and $\sigma_i = -1$ when $\beta_i \in \roots^+\setminus\roots_0^+$.
\end{itemize}
Since $1 - e_{2\sigma_i\beta_i}$ is also invertible for all $i$,
$A_{l,\nu}$ is invertible. Hence we can solve the equations
$A_{l,\nu}(m_{\Gamma_i})_i + (u_i)_i = 0$ inductively on $\nu$
and obtain a unique weak solution.

(ii) We show the statement by induction on $l$. If $l = 1$,
there is nothing to prove. For induction step, we assume that
the statement holds for $l - 1$ with $l \ge 2$. Then any weak solution
of $(\ref{eq:recursion formula})$ for all $\alpha^{\bs{i}}_m$ with
$1 \le m \le l$ is actually a solution for all $1 \le m < l$.
 Now we define $m'_{\Lambda}$ as the coefficient of $\hF_{\bs{i}}^{(\Lambda)}$ in the LHS of (\ref{eq:recursion formula}) for $\alpha^{\bs{i}}_l$. What we would like to show is $m_{\Lambda} = 0$ for all $\Lambda$. For this, it suffices to show that $(m'_{\Lambda})_{\Lambda}$ is asolution of (\ref{eq:recursion formula}) for $\alpha^{\bs{i}}_m$ with $1 \le m \le l - 1$ since the uniqueness part of (i) and $m'_{\Lambda} = 0$ for $l < \Lambda \le N - N_S$ implies $m'_{\Lambda} = 0$
for all $\Lambda$.
This follows from Proposition \ref{prop:LS relation},
which states that $\hE_{\bs{i},m}\hE_{\bs{i},l}$ can be expressed as a linear combination of $\hE_{\bs{i},l}\hE_{\bs{i},m}$ and other PBW vectors of the form $\hE_{\bs{i}}^{\Lambda^+}$ with $m < \Lambda^+ < l$.
\end{proof}
\begin{proof}[Proof of Proposition \ref{prop:uniqueness of generalized maximal vector}]
This is a special case of Lemma \ref{lemm:uniqueness and existence of partly maximal vector}, namely the case of $l = N - N_S$.
\end{proof}
\begin{defn}
Set $V = U_h^k(\frl_S)$ and $W = U_h^k(\g)$.
A generalized $S$-maximal vector $(U_{\Lambda})_{\Lambda}$ with $U_0 = 1$ is called the \emph{$S$-maximizer}.
\end{defn}

The $S$-maximizer is a universal solution of (\ref{eq:recursion formula}) in the following sense.

\begin{prop}
Let $(U_{\Lambda})_{\Lambda}$ be the $S$-maximizer and
$(m_{\Lambda})_{\Lambda}$ is a generalized $S$-maximal vector.
Then $m_{\Lambda} = U_{\Lambda}m_0$ holds for any $\Lambda$.
\end{prop}
\begin{proof}
The statements immediately follows from the definition of
generalized $S$-maximal vector and the uniqueness part of
Proposition \ref{prop:uniqueness of generalized maximal vector}.
\end{proof}

\subsection{Construction of the 2-cocycles}
Next we discuss a further property of the $S$-maximizer.
We begin with the following preparation.

\begin{lemm} \label{lemm:almost commutativity}
Let $\cl{O}$ be the subalgebra of $U_q^{\A}(\g)$ generated by
$\{\tF_{\bs{i},\alpha}\}_{\alpha \in \roots^+}$,
$\{\tE_{\bs{i},\alpha}\}_{\alpha \in \roots}$ and
$\{K_{\mu}\}_{\mu \in P}$.
\begin{enumerate}
 \item A family $\{\tF_{\bs{i}}^{\Lambda^-}K_{\mu}\tE_{\bs{i}}^{\Lambda^+}\}_{\Lambda^{\pm},\mu}$ is an $\A$-basis of $\cl{O}$. 
 \item For any $x,y \in \cl{O}$, $[x,y]_{q^{\pm 1}} \in (q - 1)\cl{O}$.
\end{enumerate}
\end{lemm}
\begin{proof}
 (i) This follows from Lemma \ref{lemm:right system case} by specializing $e_{-2\alpha}$ to $1$ for all $\alpha \in \roots^+$.

 (ii) We may assume both of $x,y$ are generators in the statement.
Moreover we may assume $x,y$ are not $K_{\mu}$ for any $\mu$ since
it is easy to see the statement in such a case.

We may assume that $x = \tE_{\bs{i},\alpha}$ with $\alpha \in \simples$ and $y = \tF_{\bs{i},\beta}$ with $\beta \in \roots^+$,
using Lemma \ref{lemm:rotation lemma}. In particular $y \in U_q^{\A}(\n^-)$.
Hence $[\aE_{\bs{i},\alpha},y]_q$
is an element of $U_q^{\A}(\frb^-)$. On the other hand $\mathrm{span}_{\A}\{\tF_{\bs{i}}^{\Lambda^-}K_{\mu}\aE_{\bs{i}}^{(\Lambda^+)}\}$ is closed under the multiplication.
Hence $[\aE_{\bs{i},\alpha},y]_q$ is in the intersection
of these subalgebras, which is spanned by
$\{\tF_{\bs{i}}^{\Lambda^-}K_{\mu}\}_{\Lambda^-,\mu}$.
This means $[x,y] \in (q_{\alpha} - q_{\alpha}^{-1})\cl{O}$ since
$x = (q_{\alpha} - q_{\alpha}^{-1})\aE_{\bs{i},\alpha}$.
\end{proof}

\begin{prop}
Let $(U_{\Lambda})_{\Lambda}$ be the $S$-maximizer.
Then $U_{\Lambda} \in h^{\abs{\Lambda}}U_h^{k}(\g\times \frl_S)$ for any $\Lambda$.
\end{prop}
\begin{proof}
We construct a weak solution of (\ref{eq:recursion formula}) for all
$\beta \in \roots^+\setminus\roots_S^+$ by another twisted parabolic induction using $U_{h,\chi}^k(\tg)$ i.e. $U_{h,\chi}^k(\tg)\tensor_{U_{h,\chi}^k(\tfrp)} V$. Note that this can be naturally embedded in
$\dqind{V}$.

In this case we consider the following equation for each $\beta \in \roots^+\setminus\roots_0^+$:
\begin{align} \label{eq:another recursion formula}
 \sum_{\Lambda} \vE_{\bs{i},\beta,(1)}\tm_{\Lambda,(W)}\tensor (\vE_{\bs{i},\beta,(2)}\vF_{\bs{i}}^{\Lambda})_{(\fru_S^-)}\tensor \pi_S((\vE_{\bs{i},\beta,(2)}\vF_{\bs{i}}^{\Lambda})_{(\frp_S)})\tm_{\Lambda,(V)}
= 0.
\end{align}
Then we can define the notion of a weak solution of (\ref{eq:another recursion formula}) in a similar manner to that of (\ref{eq:recursion formula}).
If this equation has a weak solution for any $m = \tm_0$,
we can see that the statement is valid since
we can obtain a weak solution $(m_{\Lambda})_{\Lambda}$ of (\ref{eq:recursion formula}) for all $\beta \in \roots^+\setminus\roots_S^+$
by setting $m_{\Lambda} := (q_{\alpha^{\bs{i}}} - q_{\alpha^{\bs{i}}}^{-1})^{\Lambda}[\Lambda]_{q_{\alpha^{\bs{i}}}}!\tm_{\Lambda}$. Then
Lemma \ref{lemm:uniqueness and existence of partly maximal vector} (ii) implies that this is actually a generalized $S$-maximal vector.
In particular we can see $m_{\Lambda} \in h^{\abs{\Lambda}}W\htensor V$.

The outline of the proof is almost the same as that of
the proof of Lemma \ref{lemm:uniqueness and existence of partly maximal vector} with $l = N - N_S$.
The only different point is that $A_{\nu}$, the matrix of the coefficient of $\vF_{\bs{i}}^{\Gamma'}$ in the term of $\Lambda$ with $\Pi(\Lambda\cdot\alpha^{\bs{i}}) = \Pi(\Gamma\cdot\alpha^{\bs{i}}) = \nu$,
is a lower triangular matrix with invertible diagonal entries modulo $h$. In the following, we prove this fact.

Fix $\Gamma$ and take $\Gamma'$ and $\beta = \alpha^{\bs{i}}_i$
in the same way as Proof of Lemma \ref{lemm:uniqueness and existence of partly maximal vector}.
Noting that $U_h(\g)$ is cocommutative modulo $h$, we can see that
$\Delta(\vE_{\bs{i},\beta}) = \vE_{\bs{i},\beta}\tensor 1 + K_{\beta}^{-1}\tensor \vE_{\bs{i},\beta}$ modulo $h$ and the first term
does not affect to the entries of $A_{\nu}$.
Hence it suffices to look at the coefficient of $\vF_{\bs{i}}^{\Gamma'}$ in the following:
\begin{align*}
\sum_{\Lambda} K_{\beta}^{-1}\tm_{\Lambda,(W)}\tensor [\vE_{\bs{i},\beta}, \vF_{\bs{i}}^{\Lambda}]_{q,(\fru_S^-)}\tensor \pi_S([\vE_{\bs{i},\beta}, \vF_{\bs{i}}^{\Lambda}]_{q,(\frp_S)})\tm_{\Lambda,(V)}. 
\end{align*}
We consider the following decomposition of a quantum commutator:
\begin{equation} \label{eq:expansion of quantum commutator}
\begin{aligned}
 [\vE_{\bs{i},\beta},\vF_{\bs{i}}^{\Lambda}]_q
= [\vE_{\bs{i},\beta},\vF_{\bs{i}}^{\Lambda_{<i}}]_q\vF_{\bs{i}}^{\lambda_i}\vF_{\bs{i}}^{\Lambda_{>i}}
&+ q^{(\beta,\Lambda_{<i}\cdot\alpha^{\bs{i}})}\vF_{\bs{i}}^{\Lambda_{<i}}[\vE_{\bs{i},\beta},\vF_{\bs{i},i}^{\lambda_i}]_q\vF_{\bs{i}}^{\Lambda_{>i}} \\
&+ q^{(\beta,\Lambda_{<i}\cdot\alpha^{\bs{i}} + \lambda_i\alpha^{\bs{i}}_i)}\vF_{\bs{i}}^{\Lambda_{<i}}\vF_{\bs{i},i}^{\lambda_i}[\vE_{\bs{i},\beta},\vF_{\bs{i}}^{\Lambda_{>i}}]_q
\end{aligned}
\end{equation}
Assume $\Lambda_{<i} \neq 0$. It is easy to see that the second and
third terms do not affect the coefficient of
$\vF_{\bs{i}}^{\Gamma'}$.
We show that even the first term does not modulo $h$.
By Proposition \ref{prop:mixed LS relation} and Lemma \ref{lemm:almost commutativity} we have the following expansion:
\begin{align*}
 [\vE_{\bs{i},\beta},\vF_{\bs{i}}^{\Lambda_{<i}}]_q
= \sum_{\substack{\Lambda^- < i \le \Lambda^+,\mu \\ \beta - \Lambda^+\cdot\alpha^{\bs{i}} = \Lambda_i\cdot\alpha^{\bs{i}} - \Lambda^-\cdot\alpha^{\bs{i}}}} C_{\Lambda^{\pm},\mu} \vF_{\bs{i}}^{\Lambda^-}K_{\mu}\hE_{\bs{i}}^{\Lambda^+}
\end{align*}
Then Lemma \ref{lemm:almost commutativity} allows us to
expand the first term as follows modulo $h$:
\begin{align*}
\sum_{\substack{\Lambda^- < i \le \Lambda^+,\mu \\ \beta - \Lambda^+\cdot\alpha^{\bs{i}} = \Lambda_i\cdot\alpha^{\bs{i}} - \Lambda^-\cdot\alpha^{\bs{i}}}} C_{\Lambda^{\pm},\mu} \vF_{\bs{i}}^{\Lambda^-}\vF_{\bs{i}}^{\lambda_i}\vF_{\bs{i}}^{\Lambda_{>i}}K_{\mu}\hE_{\bs{i}}^{\Lambda^+}.
\end{align*}
In this summation $\vF_{\bs{i}}^{\Gamma'}$ appears as the
$U_{h,\chi}^k({}^{\sim}\fru_S^-)$-part only when $\Lambda^- = 0$, in which case $\Lambda^+\neq 0$ since $\Lambda^+\cdot\alpha^{\bs{i}} - \Lambda^-\cdot\alpha^{\bs{i}} = \beta - \Lambda_{<i}\cdot\alpha^{\bs{i}} \neq 0$
by \cite[Theorem, p.~662]{MR1169886}. If $\Lambda^+ > N - N_S$ holds additionally,
$U_h^k(\frl_S)$-part of such a term is $0 \modulo h$ since
$\hE_{\bs{i}}^{\Lambda^+} \in hU_h^k(\frl_S)$.
If $\Lambda^- \not> N - N_S$ holds, such a term is removed by $\pi_S$.
Hence we can ignore $\Lambda$ with $\lambda_1 = \lambda_2 = \cdots = \lambda_{i - 1} = 0$.

Next assume $\Lambda_{< i} = 0$. The case of $\Lambda = \Gamma$ can be
treated in the same way as the proof of Lemma \ref{lemm:uniqueness and existence of partly maximal vector}.
When $\Lambda \neq \Gamma$, the second term does not affect
the coefficient of $\vF_{\bs{i}}^{\Gamma'}$
and the third term does only when $\lambda_i = \gamma_i - 1 < \gamma_i$.

Therefore $A_{\nu}$ is of the following form modulo $h$:
\begin{align} \label{eq:another linear equation}
\begin{pmatrix}
 [\gamma_{1,i_1}]_{q_{\beta_1}}\sigma_1(1 - e_{-2\sigma_1\beta_1})
& \cdots & 0 \\
\vdots 
& \ddots & \vdots \\
 \ast 
& \cdots & [\gamma_{l,i_{n_{\nu}}}]_{q_{\beta_{n_{\nu}}}}\sigma_{n_{\nu}}(1 - e_{-2\sigma_{n_{\nu}}\beta_{n_{\nu}}})
\end{pmatrix},
\end{align}
where $\Gamma_k = \Gamma'_k + \delta_{i_k}$ and $\beta_k = \alpha^{\bs{i}}_{i_k}$ are taken as before. Since $[\gamma_{k,i_k}]_{q_{\beta_k}}$ is invertible in our setting,
we can see the existence of a weak solution of
(\ref{eq:another recursion formula}) for all $\beta \in \roots^+\setminus\roots_S^+$.
\end{proof}
\begin{coro}
\label{coro:existence of maximal vector in tensor products}
Let $V$ be a $U_h^k(\frl_S)$-module and $W$ be a $U_h^k(\g)$-module.
Then there is a unique
$S$-maximal vector of the form $m_{(W)}\htensor (1\tensor m_{(V)}) + \cdots$ for any $m \in W\htensor V$.
\end{coro}
\begin{proof}
By the definition of the completion, we can consider
$\sum_{\Lambda} (U_{\Lambda}m)_{(W)}\tensor \hF_{\bs{i}}^{(\Lambda)}\tensor (U_{\Lambda}m)_{(V)}$, which is the $S$-maximizer.
\end{proof}
\begin{rema} \label{rema:no need of completion for finite integrable}
If $V$ and $W$ are finite integrable, there is no need to take the
completion.
\end{rema}
\begin{rema} \label{rema:maximal vectors in parabolic induction}
Since $1\tensor v$ is $S$-maximal in $\dqind V$, the uniqueness
implies that an $S$-maximal vector in $\dqind V$ must be of
the form $1 \tensor v$.
\end{rema}
\begin{defn}
Let $\rep_{h,\chi}^{\fin} \frl_S$ be the full subcategory
of $\mod{U_{h,\chi}^k(\gt)}$ consisting of objects which are isomorphic to $\dqind V$ for some $V \in \rep^{\fin}_h \frl_S$.
\end{defn}
\begin{prop} \label{prop:fundamentals on dqind}
The following facts hold:
\begin{enumerate}
 \item The functor $\dqind$ restricts to an equivalence of category
between $\rep_h^{\fin}\frl_S$ and $\rep_{h,\chi}^{\fin}\frl_S$.
 \item For $V \in \rep_h^{\fin} \frl_S$ and $W \in \rep_h^{\fin} \g$,
there is a natural isomorphism from $\dqind (W\tensor V)$ to
$W\tensor \dqind V$ such that $1 \tensor (w \tensor v)\longmapsto w\tensor (1\tensor v) + \cdots$.
\end{enumerate}
Therefore $\rep_{h,\chi}^{\fin} \frl_S$ naturally becomes a semisimple left $\rep_h^{\fin}\g$-module category.
\end{prop}
\begin{proof}
Using the universal property of $(\roots_0^+,\chi)$-parabolic induction (Lemma \ref{lemm:universal property}),
we can deduce (i) from Remark \ref{rema:maximal vectors in parabolic induction}.

Next we prove (ii). By the universal property and Corollary \ref{coro:existence of maximal vector in tensor products} we can see
that a homomorphism stated in (ii) actually exists and is unique.
Hence it suffices to show that this is an isomorphism.
Consider the following composition of homomorphisms:
\begin{align*}
W\tensor\dqind V
&\longrightarrow W\tensor \dqind(W^*\tensor W\tensor V)\\
&\longrightarrow W\tensor W^*\tensor \dqind(W\tensor V)
\longrightarrow \dqind(W\tensor V),
\end{align*}
where $k\bbh\longrightarrow W^*\tensor W$ is given by
$1\longmapsto \sum_i e^i\tensor e_i$ with $(e_i)_i$ is a basis of
$W$ and $(e^i)_i$ is its dual basis of $W^*$, and $W\tensor W^*\longrightarrow k\bbh$ is given by $w\tensor f\longmapsto f(K_{2\rho}w)$
where $\rho$ is the half sum of positive roots.
Then the composition of this homomorphism and $\dqind (W\tensor V)\longrightarrow W\tensor \dqind V$ maps $1\tensor (w\tensor v)$ to $1\tensor (w\tensor v) \modulo h$ and is isomorphism. This fact implies
the injectivity of the homomorphism in the statement.
We can see surjectivity in a similar way.
\end{proof}

We rewrite $\rep_{h,\chi}^{\fin}\frl_S$ as a $2$-cocycle twist
of $\rep_h^{\fin}\frl_S$. For $V \in \rep_h^{\fin}\frl_S$ and $W,W' \in \rep_h^{\fin}\g$, we consider the following isomorphism:
\begin{equation} \label{eq:associator for twisted induction}
\begin{aligned}
\dqind(W\tensor W'\tensor V)
&\cong W\tensor\dqind(W'\tensor V)\\
&\cong (W\tensor W')\tensor \dqind V
\cong \dqind (W\tensor W'\tensor V).
\end{aligned}
\end{equation}
Proposition \ref{prop:fundamentals on dqind} (i) implies that
this automorphism is induced from an automorphism on $W\tensor W'\tensor V$. Actually we have a stronger statement.
The centralizer of $U_h^k(\frl_S)$ embedded in $U_h^k(\g\times\g\times \frl_S)$ via $(\Delta\tensor \id)\Delta$ is denoted by
$U_h^k(\g\times\g\times\frl_S)^{\frl_S}$.

\begin{prop} \label{prop:associated 2-cocycle}
There is a unique $\Psi_{h,\chi}^k \in U_h(\g\times\g\times \frl_S)^{\frl_S}$ such that the isomorphism (\ref{eq:associator for twisted induction}) is induced from the mulitplication by $\Psi_{h,\chi}$.
Moreover $\Psi_{h,\chi}$ satisfies 
the cocyle identities
\begin{align} \label{eq:cocycle equation}
&(\id\tensor \Delta\tensor \id)(\Psi_{h,\chi})(1\tensor \Psi_{h,\chi})
= (\Delta\tensor \id\tensor \id)(\Psi_{h,\chi})(\id\tensor \id\tensor \Delta)(\Psi_{h,\chi}), \\
\label{eq:cocycle normalization}
&(\eps\tensor\id\tensor \id)(\Psi_{h,\chi})
= (\id\tensor\eps\tensor\id)(\Psi_{h,\chi}) = 1
\end{align}
and has the following expansion:
\begin{align} \label{eq:expansion of 2-cocycle}
 \Psi_{h,\chi} = 1 - h&\left\{\sum_{\beta \in \roots^+\cap\roots_0^+\setminus\roots_S^+} \frac{2d_{\beta}}{1 - \chi_{-2\beta}}E_{\bs{i},\beta}\tensor F_{\bs{i},\beta}\tensor 1 \right.\\
&\hspace{5em}\left.+ \sum_{\beta \in \roots^+\setminus\roots_0^+} \frac{2d_{\beta}\chi_{2\beta}}{\chi_{2\beta} - 1}E_{\bs{i},\beta}\tensor F_{\bs{i},\beta}\tensor 1\right\} + \cdots.
\end{align}
\end{prop}
\begin{proof}
Fix $u \in W\tensor W'\tensor V$.
Let $u'$ be the image of $1\tensor u$
at $W\tensor W'\tensor \dqind V$ in (\ref{eq:associator for twisted induction}) and $\Psi_{W,W',V} \in \mathrm{End}_{U_h^k(\frl_S)}(W\tensor W'\tensor V)$ be the isomorphism corresponding to (\ref{eq:associator for twisted induction}). Then we have
\begin{align*}
u' =
\Psi_{W,W',V}(u)_{(W\tensor W')}\tensor (1 \tensor \Psi_{W,W',V}(u)_{(V)}) + \cdots.
\end{align*}
On the other hand, by straightforward calculation, we also have
\begin{align*}
u' = \sum_{\Lambda} &((\id\tensor \Delta)(U_{\Lambda})u)_{(W)}\\
&\tensor \hF_{\bs{i}}^{(\Lambda)}(((\id\tensor \Delta)(U_{\Lambda})u)_{(W')}
\tensor (1\tensor ((\id\tensor \Delta)(U_{\Lambda})u)_{(V)}) + \cdots).
\end{align*}
Hence we can obtain $\Psi_{h,\chi}$ as follows:
\begin{align} \label{eq:explicit presentation of 2-cocycle}
 \Psi_{h,\chi} = \sum_{\Lambda} (1\tensor \hF_{\bs{i},(1)}^{\Lambda}\tensor \pi_S(\hF_{\bs{i},(2)}^{\Lambda}))(\id\tensor \Delta)(U_{\Lambda}).
\end{align}
Note that this summation is well-defined since $U_{\Lambda} \in h^{\abs{\Lambda}}U_h^k(\g\times\g\times \frl_S)$.

The cocycle identities (\ref{eq:cocycle equation}), (\ref{eq:cocycle normalization}) can be obtained as consequences of direct calculations.
To prove the expansion,
we use the formula (\ref{eq:explicit presentation of 2-cocycle}).
Since $U_{\Lambda} \in h^{\abs{\Lambda}}U_h^k(\g\times\g\times\frl_S)$
and
\begin{align*}
 \Delta(\hF_{\bs{i},\beta}) = 
\begin{cases}
 F_{\bs{i},\beta}\tensor 1 + K_{\beta}^{-1}\tensor F_{\bs{i},\beta} & \modulo h \quad(\beta \in \roots^+\cap\roots_0^+), \\
 \chi_{2\beta}F_{\bs{i},\beta}\tensor 1 + K_{\beta}^{-1}\tensor F_{\bs{i},\beta} & \modulo h \quad(\beta \in \roots^+\setminus \roots_0^+),
\end{cases}
\end{align*}
it suffices to calculate the following summation:
\begin{align*}
 \sum_{\beta \in \roots^+\cap\roots_0^+} (1\tensor F_{\bs{i},\beta}\tensor 1)(\id\tensor \Delta)(U_{\beta}) +  \sum_{\beta \in \roots^+\setminus\roots_0^+} \chi_{2\beta}(1\tensor F_{\bs{i},\beta}\tensor 1)(\id\tensor \Delta)(U_{\beta}),
\end{align*}
where $U_{\beta} := U_{\delta_k}$ with $\beta = \alpha^{\bs{i}}_k$.
To determine $U_{\beta}$ modulo $h$
for each $\beta = \alpha^{\bs{i}}_k$,
we go back to (\ref{eq:another linear equation}) and show that
entries in the $i$-th row of $A_{\Pi(\beta)}$, where $\Gamma_i = \delta_k$,
are $0$ modulo $h$ except the diagonal entry.
Let $\Lambda$ be a multi-index
with $\Pi(\Lambda\cdot\alpha^{\bs{i}}) = \Pi(\beta)$
and $\Lambda \ge k$ and $\Lambda \neq \delta_{i_k}$.
Then we have $\Lambda > k$, in which case we can use Proposition \ref{prop:mixed LS relation} to obtain the following expansion:
\begin{align*}
 [\vE_{\bs{i},k}, \vF_{\bs{i}}^{\Lambda}]
= \sum_{\Lambda^+ < k \le \Lambda^-} C_{\Lambda^{\pm}} \vF_{\bs{i}}^{\Lambda^-}\vE_{\bs{i}}^{(\Lambda^+)}.
\end{align*}
Since $\pi_S(\vE_{\bs{i}}^{\Lambda^+}) = 0$ whenever $\Lambda^+ \neq 0$,
it suffices to look at the terms with $\Lambda^+$. Moreover $\Gamma_i = 0 + \delta_k$ implies that we may additionally assume $\Pi(\Lambda^-\cdot\alpha^{\bs{i}}) = 0$.
Such a term is $0$ modulo $h$ in $U_h^k(\frp_S)$ if $\Lambda^- \neq 0$, which is automatically satisfied since
$\Lambda\cdot\alpha^{\bs{i}} \neq \alpha^{\bs{i}}_k$
by \cite[Theorem, p.662]{MR1169886} and \cite[Chapter VI, Proposition 19]{MR1890629}.


Therefore we have
\begin{align*}
U_{\beta} 
= -\sigma\frac{q_{\beta} - q_{\beta}^{-1}}{1 - \chi_{-2\sigma\beta}}E_{\bs{i},\beta}K_{\beta}^{-1}\tensor 1
= -\sigma\frac{2d_{\beta}}{1 - \chi_{-2\sigma\beta}}E_{\bs{i},\beta}\tensor 1 \quad\modulo h
\end{align*}
where $\sigma = 1 $ when $\beta \in \roots_0^+$ and $\sigma = -1$ when
$\beta \not\in \roots_0^+$ and obtain the expansion (\ref{eq:expansion of 2-cocycle}).
\end{proof}

\subsection{Construction of quantum semisimple coadjoint orbits}

Let $G$ be a $1$-connected complex Lie group with the Lie algebra $\g$
and $L_S$ be its connected complex Lie subgroup with the Lie algebra $\frl_S$.
In this subsection we use $\Psi_{h,\chi}$ to construct
a quantum semisimple coadjoint orbit, which is a deformation quantization of $\cl{O}(L_S \backslash G)$ with an action of $U_h(\g)$.
Let $\irr_h\g$ be the set of equivalence classes of
irreducible finite integrable
$U_h^k(\g)$ modules. We also fix a representative $\pi$ of each class
and sometimes identify $\irr_h \g$ as the set of all representatives.

For a $U_h^k(\g)$-module $V$,
its $U_h^k(\frl_S)$-invariant part is denoted by $V^{\frl_S}$.

Set $\cl{O}_{h,\chi}^k(L_S\backslash G)$ as follows:
\begin{align} \label{eq:spectral decomposition}
 \cl{O}_{h,\chi}^k(L_S\backslash G) = \bigoplus_{\pi \in \irr_h \g}
\Hom_{U_h^k(\frl_S)}(V_{\pi},k\bbh)\tensor V_{\pi}.
\end{align}
Note that there is a unique linear map $\map{p_V}{\Hom_{U_h^k(\frl_S)}(V,k\bbh)\tensor V}{\cl{O}_{h,\chi}^k(L_S\backslash G)}$ for each finite integrable module $V$
which is natural in the 
sense that $p_V(f\circ T\tensor v) = p_W(f\tensor T(v))$ 
for any $f \in \Hom_{U_h^k}(\frl_S,k\bbh),v \in V$ and $T \in \Hom_{\g_h}(V,W)$.
For ease to read, we use $f\tensor v$ to represent
$p_V(f\tensor v)$.

Using this we construct a product $\ast_{h,\chi}$
on $\cl{O}_{h,\chi}^k(L_S\backslash G)$ defined as follows by using
$\Psi_{h,\chi,0} = (\id\tensor\id\tensor\eps)(\Psi_{h,\chi})$,
which is thought as a $U_h^k(\frl_S)$-homomorphism:
\begin{align*}
 (f\tensor v)\ast_{h,\chi}(g\tensor w) = (f\tensor g)\Psi_{h,\chi,0}^{-1}\tensor (v\tensor w)
\end{align*}

This product is actually associative. To see this
consider the product of three elements in different order:
\begin{align*}
 ((f\tensor u)\ast_{h,\chi}(g\tensor v))&\ast_{h,\chi}(h\tensor w) \\
&= (f\tensor g\tensor h)(\Psi_{h,\chi,0}^{-1}\tensor 1)(\Delta\tensor \id)(\Psi_{h,\chi,0}^{-1})\tensor (u\tensor v\tensor w),\\
 (f\tensor u)\ast_{h,\chi}((g\tensor v)&\ast_{h,\chi}(h\tensor w))\\
&= (f\tensor g\tensor h)(1\tensor \Psi_{h,\chi,0}^{-1})(\id\tensor \Delta)(\Psi_{h,\chi,0}^{-1})\tensor (u\tensor v\tensor w).
\end{align*}
Since $hx = \eps(x)h$ for all $x \in U_h^k(\frl_S)$, we can see that
$(f\tensor g \tensor h)(\Psi_{h,\chi,0}^{-1}\tensor 1) = (f\tensor g\tensor h)\Psi_{h,\chi}^{-1}$. Then the associativity follows from the cocycle identity (\ref{eq:cocycle equation}).

When $k = \C$, $\cl{O}_{h,\chi}(L_S\backslash G)$ has a natural structure of deformation quantizations of $\cl{O}(L_S\backslash G)$.
To determine the semi-classical limit of $\cl{O}_{h,\chi}(L_S \backslash G)$,
we recall the identification of $U_h(\g)$ and the twist of
$U(\g)\bbh$.
For any subset $S \subset \simples$ we can take
an isomorphism $\map{T}{U(\g)\bbh}{U_h(\g)}$ as topological $\C\bbh$-algebra such that
\begin{enumerate}
 \item $T(U(\frl_S)\bbh) = U_h(\frl_S)$.
 \item There is $F \in U(\g\times \g)\bbh$ such that $F = 1 \modulo h$
and $(T\tensor T)^{-1}\circ\Delta\circ T(x) = F\Delta(x)F^{-1}$
\end{enumerate}
The standard $r$-matrix is given by the following formula:
\begin{align*}
 r = \sum_{\alpha \in \roots^+} d_{\alpha}E_{\alpha}\wedge F_{\alpha},
\end{align*}
here $\g\wedge \g$ is embedded into $\g\tensor \g$ by
$x\wedge y \longmapsto \frac{1}{2}(x\tensor y - y\tensor x)$.
Then the anti-symmetric part of the first coefficient of $F$
is given by $-r$.

Under this identification of $U_h(\g)$ with $U(\g)\bbh$,
$V$ is finite integrable $U_h(\g)$-module if and only if
it is isomorphic to $L\tensor \C\bbh$ where $L$ is
finite dimensional representation of $\g$. In particular we can
identify $\irr_h\g$ with $\irr\g$.

Using this picture of $U_h(\g)$ we can compare
$\cl{O}_{h,\chi}(L_S\backslash G)$ with $\cl{O}(L_S\backslash G)\bbh$.
At first one should note that we can decompose $\cl{O}(L_S\backslash G)$ as follows:
\begin{align*}
&\cl{O}(L_S\backslash G) = \bigoplus_{[\pi] \in \irr \g} \Hom_{\frl_S}(L_{\pi},\C)\tensor L_{\pi}, \\
&(f\tensor v)(g\tensor w) = (f\tensor g)\tensor (v\tensor w).
\end{align*}
Here the isomorphism is given by $f\tensor v \longmapsto f(\pi(\tend)v)$.

Hence we have a canonical identification of $\cl{O}_{h,\chi}(L_S\backslash G)$ with $\cl{O}(L_S\backslash G)\bbh$. To interpret $\ast_{h,\chi}$ as a product on $\cl{O}(L_S\backslash G)\bbh$,
one should note that $\Hom_{\g}(W,W'\tensor W'')\cong \Hom_{U_h(\g)}(W,W'\tensor W'')$ via $T \longmapsto F\circ T$. Therefore we have
\begin{align*}
 (f\tensor v)\ast_{h,\chi}(g\tensor w) = (f\tensor g)\Psi_{h,\chi,0}^{-1}F\tensor F^{-1}(v\tensor w).
\end{align*}

Now we can obtain the semi-classical limit of $\cl{O}_{h,\chi}(L_S\backslash G)$ by a direct calculation. For $\theta \in \g\wedge \g$, the right (resp. left) invariant bi-vector field corresponding to $\theta$ is denoted by $L_*(\theta)$ (resp. $R_*(\theta)$).

\begin{prop}
The algebra $\cl{O}_{h,\chi}(L_S \backslash G)$ is an equivariant deformation quantization of $\cl{O}(L_S \backslash G)$. The associated Poisson structure
is given by
\begin{align*}
 \{F,G\} = (L_*(\pi) + R_*(r))(F,G)
\end{align*}
where
\begin{align*}
 \pi = \sum_{\beta \in \roots^+\cap\roots_0^+} d_{\beta}\frac{1 + \chi_{-2\beta}}{1 - \chi_{-2\beta}}E_{\beta}\wedge F_{\beta} + \sum_{\beta \in \roots^+\setminus\roots_0^+} d_{\beta}\frac{\chi_{2\beta} + 1}{\chi_{2\beta} - 1}E_{\beta}\wedge F_{\beta}.
\end{align*}
\end{prop}
\section{Universal families of $2$-cocycles}
\label{section:globalization of the 2-cocycle}

\subsection{The toric variety associated to a root system}
For general treatment of toric varieties, see \cite{MR1234037,telen2022}.
Let us recall the definition of
the fan associated to a root system. For
a root system $\roots$, $Q^{\vee}$ denotes the dual lattice $\Hom_{\Z}(Q,\Z)$. We also set $\h_{\R}^*$ (resp. $\h_{\R}^{**}$) as
$Q\tensor_{\Z}\R$ (resp. $Q^{\vee}\tensor_{\Z}\R$), which are
thought as an $\R$-linear subspace of $\h^*$ (resp. $\h^{**}$).

Then the fan associated to $\roots$ is the collection
$\Sigma_{\roots}$ of the following cones:
\begin{align*}
 \sigma_{\roots_0^+,S}
:= \{\lambda \in \h^*_{\R}\mid \lambda|_{\roots_0^+} \ge 0, \lambda|_{S_0} = 0 \},
\end{align*}
where $\roots_0^+$ is a positive system and $S_0$ is a set of its simple roots. It is not difficult to see that this collection is actually a
fan i.e. satisfies the following conditions:
\begin{enumerate}
 \item For any $\sigma \in \Sigma_{\roots}$, $\sigma \cap -\sigma = \{0\}$ holds.
 \item For any $\sigma \in \Sigma_{\roots}$ and $m \in \sigma^{\vee} = \{x \in \h_{\R}\mid \lambda(x) \ge 0 \text{ for } \lambda \in \sigma\}$, the face $\sigma_m := \{\lambda \in \sigma \mid \lambda(m) = 0\}$ is also contained in $\Sigma_{\roots}$.
 \item For any $\sigma, \sigma' \in \Sigma_{\roots}$, $\sigma \cap \sigma'$ is of the form $\sigma_m$ for some $m \in \sigma^{\vee}$.
\end{enumerate}
Let $k$ be a commutative ring. The affine toric scheme $U_{\sigma}(k)$ associated to $\sigma \in \Sigma_{\roots}$ is defined as
$\spec k[\sigma^{\vee}\cap 2Q]$, here we use $2Q$ instead of $Q$ by a notational reason. Since any face of $\sigma$
is of the form $\sigma_m$ with $m \in \sigma^{\vee}\cap 2Q$,
a face $\tau$ of $\sigma$ gives a canonical open subscheme, namely
$\spec k[\sigma^{\vee}\cap 2Q]_{e_m} = \spec k[\tau^{\vee}\cap 2Q]$ for $\tau = \sigma_m$ with $m \in \sigma^{\vee}\cap 2Q$.
Then the toric scheme associated to $\Sigma_{\roots}$ over $k$,
denoted by $X_{\roots}(k)$, can be obtained by gluing $\{U_{\sigma}\}_{\sigma \in \Sigma_{\roots}}$ along with the gluing data
$\{U_{\sigma\cap \tau} \rightarrow U_{\sigma}\}_{\sigma,\tau \in \Sigma_{\roots}}$. 

\subsection{The moduli space of equivariant Poisson brackets on $L_S\backslash G$}
\label{subsection:the moduli space}
Let $\cl{O}(X_{L_S\backslash G})$ be the quotient of $\C[\phi_{\alpha}]_{\alpha \in \roots\setminus\roots_S}$ divided by the following relations:
\begin{enumerate}
 \item $\phi_{-\alpha} = -\phi_{\alpha}$ for $\alpha \in \roots\setminus\roots_S$, \label{relation:(i)}
 \item $\phi_{\alpha}\phi_{\beta} + 1 = \phi_{\alpha + \beta}(\phi_{\alpha} + \phi_{\beta})$ when $\alpha,\beta,\alpha + \beta \in \roots\setminus\roots_S$, \label{relation:(ii)}
 \item $\phi_{\alpha} = \phi_{\beta}$ when $\alpha, \beta \in \roots\setminus\roots_S$ and $\alpha - \beta \in \roots_S$. \label{relation:(iii)}
\end{enumerate}
As pointed out in \cite{MR1817512, MR2141466, MR2349621},
this is the moduli of Poisson structures on $L_S\backslash G$
equivariant with respect to the Poisson-Lie group structure on $G$.
Actually we can associate an equivariant Poisson bracket
$\{\tend,\tend\}_{\phi}$ to each $\phi \in X_{L_S\backslash G}$ as follows:
\begin{align*}
 \{F,G\}_{\phi} = (L_*(\pi_{\phi}) + R_*(r))(F,G),
\text{ where } \pi_{\phi}
= \sum_{\beta \in \roots^+\setminus\roots_S^+} d_{\beta}\phi_{\beta}E_{\beta}\wedge F_{\beta}.
\end{align*}

In the following we describe a canonical immersion from $X_{L_S\backslash G}(k)$ into the toric scheme $X_R(k)$ for an arbitrary commutative ring $k$.

Let $\cl{O}(X_{L_S\backslash G}(k))$ be the quotient of $k[\psi_{\alpha}]_{\alpha \in \roots\setminus\roots_S}$ divided by the following relations:
\begin{enumerate}
\renewcommand{\labelenumi}{(\roman{enumi}${}_{\psi}$)}
 \item $\psi_{\alpha} + \psi_{-\alpha} = 1$ for $\alpha \in \roots\setminus\roots_S$, \label{relation:(i)'}
 \item $\psi_{\alpha}\psi_{\beta} = \psi_{\alpha + \beta}(\psi_{\alpha} + \psi_{\beta} - 1)$ when $\alpha,\beta,\alpha + \beta \in \roots\setminus\roots_S$, \label{relation:(ii)'}
 \item $\psi_{\alpha} = \psi_{\beta}$ when $\alpha ,\beta \in \roots\setminus\roots_S$ satisfy $\alpha - \beta \in \aspan_{\Z} S$. \label{relation:(iii)'}
\end{enumerate}
Then we set $X_{L_S\backslash G}(k) = \spec \cl{O}(X_{L_S\backslash G}(k))$. Note that these relations coincide with
(\ref{relation:(i)}), (\ref{relation:(ii)}), (\ref{relation:(iii)})
under the transformation $\phi_{\alpha} = 2\psi_{\alpha}-1$ if $2 \in k^{\times}$.
In particular we can canonically identify $X_{L_S\backslash G}(\C)$
with $X_{L_S\backslash G}$.
\begin{lemm}
 For a positive system $\roots_0^+$,
we set $f_{\roots_0^+} = \prod_{\alpha \in\roots_0^+\setminus\roots_S}\psi_{\alpha}$ and $U(\roots_0^+) := \{\frp \in X_{L_S\backslash G}(k)\mid f_{\roots_0^+} \not\in \frp\}$. Then $\{U(\roots_0^+)\}_{\roots_0^+ \supset \roots_S^+}$
is an open covering of $X_{L_S\backslash G}(k)$.
\end{lemm}
\begin{proof}
Take a prime ideal $\frp$. What we have to show is the existence of a
positive system $\roots_0^+$ containing $\roots_S^+$ such that
$\psi_{\alpha} \not\in \frp$ for all $\alpha \in \roots_0^+\setminus \roots_S^+$,
which is equivalent to $\psi_{\alpha} \neq 0$ in $\cl{O}(X_{L_S \backslash G}(k))/\frp$. Since this algebra is
a domain, the following lemma can be applied.
\end{proof}
\begin{lemm}[c.f. {\cite[Lemma 11]{MR2141466}}]
 Let $k$ be a domain and consider 
$\psi = \{\psi_{\alpha}\}_{\alpha} \in k^{\roots\setminus\roots_S}$
satisfying (\ref{relation:(i)'}${}_{\psi}$), (\ref{relation:(ii)'}${}_{\psi}$), (\ref{relation:(iii)'}${}_{\psi}$).
Then there is a positive system $\roots_0^+$ contained in
$\roots_S^+ \cup \{\alpha \in \roots\setminus\roots_S\mid \psi_{\alpha}\neq 0\}$.
\end{lemm}
\begin{proof}
At first we show that $P :=\roots_S \cup \{\alpha \in \roots\setminus\roots_S\mid \psi_{\alpha}\neq 0\}$ contains a positive system.
By \cite[Chapter VI, Proposition 20]{MR1890629},
it suffices to show that the following properties hold:
\begin{enumerate}
\renewcommand{\labelenumi}{(\alph{enumi})}
 \item For any $\alpha \in \roots$, either of $\alpha \in P$ or $-\alpha \in P$ holds.
 \item For any $\alpha,\beta \in P$, $\alpha + \beta \in P$ whenever
$\alpha + \beta \in \roots$.
\end{enumerate}
The condition (a) follows from (\ref{relation:(i)'}${}_{\psi}$). To see (b),
we assume either of $\alpha$ or $\beta$ is in $\roots_S$ at first. 
Then (\ref{relation:(iii)'}${}_{\psi}$) shows (b) in this case. Next we assume that both of them are not in $\roots_S$ and $\alpha + \beta \in \roots$.
Moreover we also assume $\alpha + \beta \not\in \roots_S$ since there
is nothing to prove when $\alpha + \beta \in \roots_S$. In this case,
using (\ref{relation:(ii)'}${}_{\psi}$), we have
\[
 \psi_{\alpha + \beta}(\psi_{\alpha} + \psi_{\beta} - 1) = \psi_{\alpha}\psi_{\beta}\neq 0
\]
since $k$ is a domain. Hence we have
$\psi_{\alpha + \beta}\neq 0$, which implies $\alpha + \beta \in P$.

Take a positive system $\roots_1^+$ of $\roots$ contained in $P$.
Then $\roots_1^+\cap \roots_S$ is also a positive system of $\roots_S$. Hence there is $w$, a product of $s_{\alpha}$ with $\alpha \in \roots_S$, such that $w(\roots_1^+\cap \roots_S) = \roots^+\cap \roots_S$. Now $\roots_0^+ := w(\roots_1^+)$ meets the 
required conditions since $w$ preserves a decompostion 
$\roots = \roots_S\sqcup (\roots\setminus\roots_S)$ and
$\psi_{w(\alpha)} = \psi_{\alpha}$ holds for $\alpha \in \roots\setminus\roots_S$.
\end{proof}
\begin{lemm} \label{lemm:local isomorphism}
Let $\cl{O}(U(\roots_0^+))$ be the localization of
$\cl{O}(X_{L_S \backslash G}(k))$ by $f_{\roots_0^+}$
and
$\cl{B}_{\roots_0^+}$ be the localization of $k[2Q_0^-]/(e_{-2\alpha} - 1)_{\alpha \in \roots_S^+}$ by $(1 - e_{-2\alpha})_{\alpha \in \roots_0^+\setminus\roots_S^+}$.
Then there is a canonical isomorphism between these algebras,
as constructed in the proof below.
\end{lemm}
\begin{proof}
We give the constructions of $\map{u}{\cl{O}(U(\roots_0^+))}{\cl{B}_{\roots_0^+}}$ and $\map{v}{\cl{B}_{\roots_0^+}}{\cl{O}_{X_{L_S \backslash G}}(U(\roots_0^+))}$ which are mutually inverse. It suffices
to give the images of the generators $\{\psi_{\alpha}\}_{\alpha \in \roots_0^+\setminus\roots_S^+}$ and $\{e_{-2\alpha}\}_{\alpha \in \roots_0^+\setminus\roots_S^+}$:
\begin{align} \label{eq:transformation of variables}
 u(\psi_{\alpha}) = \frac{1}{1 - e_{-2\alpha}}, \quad
 v(e_{-2\alpha}) = 1 - \frac{1}{\psi_{\alpha}}.
\end{align}
It can be checked by direct calculations that these elements
satisfy the required relations.
\end{proof}

Now we have a locally closed immersion $X_{L_S\backslash G}(k) \subset X_{\roots}(k)$.

\begin{prop} \label{prop:embedding of the moduli}
There is a canonical immersion from $X_{L_S\backslash G}(k)$ to $X_{\roots}(k)$, as constructed in the proof below.
\end{prop}
\begin{proof}
We can obtain a closed subscheme of $X_{\roots}(k)$ by
gluing $\spec k[\sigma^{\vee}\cap 2Q]/(e_{-2\alpha} - 1)_{\alpha \in \sigma^{\vee}\cap Q_S}$ and also obtain its open subscheme by gluing
$\spec \cl{B}_{\roots_0^+}$ over all $\roots_0 \supset \roots_S^+$.
Hence it suffices to check the compatibility of the isomorphisms
constructed in Lemma \ref{lemm:local isomorphism}. Since
$\spec \cl{B}_{\roots_0^+} \cap \spec \cl{B}_{\roots_1^+} = \spec \cl{B}_{\roots_0^+,\roots_1^+}$, where
$\cl{B}_{\roots_0^+,\roots_1^+}$ is the localization of $k[2Q_0^- + 2Q_1^-]/(1 - e_{-2\alpha})_{\alpha \in \roots_S^+}$ by $(1 - e_{-2\alpha})_{\alpha \in \roots_0^+\cup\roots_1^+\setminus \roots_S^+}$,
it suffices to construct an isomorphism between $\cl{B}_{\roots_0^+,\roots_1^+}$ and $\cl{O}(X_{L_S\backslash G}(k))_{f_{\roots_0^+}f_{\roots_1^+}}$ which makes the following diagram commutative:
\begin{align*}
 \xymatrix{
\cl{O}(X_{L_S\backslash G}(k))_{f_{\roots_0^+}} \ar[r] \ar[d]& \cl{O}(X_{L_S\backslash G}(k))_{f_{\roots_0^+}f_{\roots_1^+}} \ar[d]& \cl{O}(X_{L_S\backslash G}(k))_{f_{\roots_1^+}} \ar[l] \ar[d] \\
\cl{B}_{\roots_0^+} \ar[r] & \cl{B}_{\roots_0^+,\roots_1^+} & \cl{B}_{\roots_1^+} \ar[l].
}
\end{align*}
This can be proven in the completely same way as the proof of Lemma \ref{lemm:local isomorphism}.
\end{proof}

\subsection{The sheaf of deformed quantum enveloping algebra on $X_{\roots}(k)$}

In this subsection we construct a sheaf of algebras on the toric
scheme using deformed quantum enveloping algebras. 
For this purpose we have to compare deformed quantum envelpoing algebras based on different positive systems.

We begin with some preparations on root systems.
We say that $L \subset \roots^+$ is admissible if it is of the form
$\roots^+\setminus\roots_0^+$ for some positive system $\roots_0^+$,
or equivalently of the form $\roots^+\cap \roots_1^-$ for some positive system $\roots_1^+$.
The admissibility can be characterized by
the following conditions (\cite[Theorem, p.~663]{MR1169886}):
\begin{enumerate}
 \item $\alpha + \beta \in L$ for $\alpha, \beta \in L$.
 \item If $\alpha, \beta \in \roots^+$ and $\alpha + \beta \in L$,
either of $\alpha \in L$ or $\beta \in L$ holds.
\end{enumerate}
Note that $\roots^+\setminus L$ is also admissible when $L$ is
admissible and that $L$ contains a simple root if $L \neq \emptyset$.

\begin{lemm} \label{lemm:span of admissible}
Let $L \subset \roots^+$ be an admissible subset. Then there is
a unique subset $S \subset \simples$ such that
$\aspan_{\Z} L \cap \roots^+ = \roots_S^+$.
\end{lemm}
In the following proof, $\Supp \alpha := \{\eps \in \simples\mid c_{\eps} \neq 0\}$ where $\alpha = \sum_{\eps \in \simples} c_{\eps}\eps \in Q$.
\begin{proof}
By induction on $\abs{L}$. If $\abs{L} = 0, 1$, there is nothing to
prove since $L = \empty, \{\eps\}$ for some $\eps \in \simples$ respectively.
Next we proceed to the induction step. At first we show the
following claim: if $A \subset \roots^+$ and $S \subset \simples$
satisfy $\aspan_{\Z} A \cap \roots^+ = \roots_S^+$,
$\aspan_{\Z} A \cup \{\eps\} = \roots_{S\cup\{\eps\}}^+$ also holds
for any $\eps \in \simples$. This can be checked easily if $\eps \in S$. Hence we may assume $\eps \in \simples\setminus S$. Since $S$ is
a union of $\{\Supp \alpha\}_{\alpha \in A}$, it is easy
to see that $\aspan_{\Z} (A\cup\{\eps\}) \cap \roots^+ \subset \roots_{S\cup\{\eps\}}^+$. The other inclusion is trivial.

Now we take a simple root $\eps \in L$. Then $L\setminus\{\eps\}$
is admissible in $s_{\eps}(\roots^+)$. Hence there is a subset $S_0 \subset \simples$ with $\aspan_{\Z} (L\setminus\{\eps\}) \cap s_{\eps}(\roots^+) = s_{\eps}(\roots^+)_{s_{\eps}(S_0)}$. Since $-\eps$ is simple
in $s_{\eps}(\roots^+)$, the claim above implies that
$\aspan_{\Z} L \cap s_{\eps}(\roots^+) = s_{\eps}(\roots^+)_{s_{\eps}(S_0)\cup \{-\eps\}}$. Since $\aspan_{\Z} L$ is invariant under $s_{\eps}$, we can obtain the statement for $L$ with $S = S_0 \cup \{\eps\}$.
\end{proof}

\begin{coro} \label{coro:span of admissible}
Let $\roots_0^+$ and $\roots_1^+$ be positive systems of $\roots$.
Then there is $m \in \h_{\R}^{**}$ with the following conditions:
\begin{enumerate}
 \item $\{\alpha \in \roots\mid m(\alpha) > 0\} \subset \roots_0^+\cap \roots_1^+$,
 \item $\{\alpha \in \roots\mid m(\alpha) = 0\} = \aspan_{\Z}\roots_0^+\triangle \roots_1^+$,
 \item $\{\alpha \in \roots\mid m(\alpha) < 0\} \subset \roots_0^-\cap \roots_1^-$.
\end{enumerate}
\end{coro}
\begin{proof}
We may assume $\roots_1^+ = \roots^+$. Take $S \subset \simples$
by using Lemma \ref{lemm:span of admissible} for $L = \roots^+\setminus\roots_0^+$. Then $m = \sum_{\eps \in \simples\setminus S} \eps^{\vee}$ satisfies the condition.
\end{proof}

\begin{lemm} \label{lemm:compatible reduced expressions}
Let $\roots_0^+$ and $\roots_1^+$ be positive systems of $\roots$.
There are reduced expressioins $s_{\bs{i}}$ and $s_{\bs{j}}$
of $w_0$, compatible with $\roots_0^+$ and $\roots_1^+$ respectively,
and positive integers $1 \le k < l \le N$ with the following property:
\begin{enumerate}
 \item For $1 \le n \le k$, $\alpha^{\bs{i}}_n = \alpha^{\bs{j}}_n \in \roots^+\setminus(\roots_0^+\cup\roots_1^+)$.
 \item For $k < n \le l$, $\alpha^{\bs{i}}_n,\alpha^{\bs{j}}_n \in \aspan_{\Z} \roots_0^+\triangle\roots_1^+$.
 \item For $l < n \le N$, $\alpha^{\bs{i}}_n = \alpha^{\bs{j}}_n \in \roots^+\cap \roots_0^+\cap \roots_1^+$.
\end{enumerate}
\end{lemm}
\begin{proof}
By using $m \in \h_{\R}^{**}$ in Corollary \ref{coro:span of admissible}, we can see the admissiblilty of the following sets.
\begin{align*}
 L^{\pm} := \roots^+\cap \roots_0^{\pm}\cap \roots_1^{\pm}\setminus V = \{\alpha \in \roots^+\mid \pm m(\alpha) > 0\}.
\end{align*}
Take $w^{\pm} \in W$ such that $L^{\pm} = \roots^+\setminus w^{\pm}(\roots^+)$ and also take $w,w' \in W$ such that $\roots_0^+ = w(\roots^+)$ and $\roots_1^+ = w'(\roots^+)$.
Consider a reduced expression of $w^-$. This
expression can be extended to reduced expressions $w$ and
$w'$. Similarly a reduced expression of $w^+$ can be extended to
reduced expressions of $ww_0$ and $w'w_0$. Then we obtain enumerations $\{\alpha^i_n\}$ on $\roots^+\setminus\roots_i^+$ and $\{\beta^i_m\}$ on $\roots^+\cap \roots_i^+$ from these reduced expressions. 
Then, as a consequence of \cite[Theorem, p.662]{MR1169886},
we can find reduced expressions $s_{\bs{i}}$ and $s_{\bs{j}}$ which
satisfy the following conditions respectively:
\begin{align*}
 \alpha^{\bs{i}}_n =
\begin{cases}
 \alpha^0_n & (1 \le n \le \ell(w)),\\
 \beta^0_{N - n + 1} & (\ell(w) < n \le N),
\end{cases}
\quad
 \alpha^{\bs{j}}_n =
\begin{cases}
 \alpha^1_n & (1 \le n \le \ell(w')), \\
 \beta^1_{N - n + 1} & (\ell(w') < n \le N).
\end{cases}
\end{align*}
By construction these reduced expressions satisfy the required conditions.
\end{proof}
\begin{rema}
 Let $w^-$ be as defined in the proof above. Then $(w^-)^{-1}(\roots^+\cap \aspan_{\Z} \roots_0^+\triangle\roots_1^+)$ is of the form $\roots_S^+$ for some $S \subset \simples$ since the image is $\{\alpha \in \roots^+\mid m(w^-(\alpha)) = 0\}$ and its complement in $\roots^+$
is $\{\alpha \in \roots^+\mid m(w^-(\alpha)) > 0\}$.
This fact implies that the linear span
of $\{F_{\bs{i}}^{(\Lambda)}\}_{k <\Lambda \le l}$
coincides with the linear span of $\{F_{\bs{j}}^{(\Lambda)}\}_{k < \Lambda \le l}$.
\end{rema}

\begin{prop} \label{prop:comparison of different dQEAs}
Let $k$ be an $\cl{A}$-algebra and $\map{\chi}{2Q_0^- + 2Q_1^-}{k}$
be a character. Then there is a canonical isomorphism between
$U_{q,\chi_0}^k(\gt)$ and $U_{q,\chi_1}^k(\gt)$ as constructed in
the proof, where $\chi_i$ is the restriction of $\chi$ on $2Q_i^-$.
\end{prop}
\begin{proof}
Set $U_i = U_{q,e_i}^{\A[2Q_0^- + 2Q_1^+]}(\gt)$.
It suffices to show $U_0 = U_1$ in $U_q^{\A[P]}(\g)$,
where $\map{e}{2Q_0- + 2Q_1}{\A[2Q_0^- + 2Q_1^-]}$ is the
canonical character.

Take a reduced expressions $s_{\bs{i}}$ and $s_{\bs{j}}$ of $w_0$
as Lemma \ref{lemm:compatible reduced expressions}. Since
$e_{2\alpha}$ with $\alpha \in \roots^+\cap \aspan_{\Z} \roots_0^+\triangle\roots_1^+$ is invertible in $\A[2Q_0^- + 2Q_1^-]$,
$U_0$ is generated by the
following elements:
\begin{align*}
\{(q_{\alpha} - q_{\alpha}^{-1})\aE_{\bs{i},\alpha}\}_{\alpha \in\roots^+},\quad
\{e_{-\lambda}\aK_{\lambda}\}_{\lambda \in P},\quad
\begin{cases}
 e_{2\alpha^{\bs{i}}_n}\aF_{\bs{i},n}^{(r)} & (1 \le n \le k), \\
 \aF_{\bs{i},n}^{(r)} & (k < n \le l), \\
 \aF_{\bs{i},n}^{(r)} & (l < n \le N).
\end{cases}
\end{align*}
The same thing can be said for $U_1$. By the property of $s_{\bs{i}}$
and $s_{\bs{j}}$, only different generators are $\aF_{\bs{i},n}^{(r)}$
and $\aF_{\bs{j},n}^{(r)}$ with $k < n \le l$. But the algebras
generated by these elements are independent of $\bs{i}$ and $\bs{j}$ as stated in the remark above.
Hence we have $U_0 = U_1$.
\end{proof}
\begin{rema}
Take another positive system $\roots_2^+$ and consider a character
$\chi$ from $2Q_0^- + 2Q_1^- + 2Q_2^-$. In this case we have
three isomorphisms $U_{q,\chi_0}^k(\gt) \cong U_{q,\chi_1}^k(\gt)$,
$U_{q,\chi_1}^k(\gt) \cong U_{q,\chi_2}^k(\gt)$ and $U_{q,\chi_0}^k(\gt) \cong U_{q,\chi_2}^k(\gt)$.
By construction of the isomorphisms, we can see the commutativity of
the following diagram:
\begin{align*}
 \xymatrix
{
U_{q,\chi_0}^k(\gt)\ar[rr]^-{\cong} \ar[rd]_-{\cong} && U_{q,\chi_1}^k(\gt) \ar[ld]^-{\cong}\\
& U_{q,\chi_2}^k(\gt) &
}
\end{align*}
\end{rema}

Now we can define a sheaf of deformed quantum enveloping algebras on $X_{\roots}(k)$ and $X_{L_S\backslash G}(k)$.
\begin{defn}
Let $k$ be an $\A$-algebra.
We define $\rsfs{U}_{q,X_{\roots}}^k(\gt)$,
a sheaf of algebras on $X_{\roots}(k)$,
by gluing $\{U_{q,e}^{k[2Q_0^-]}(\gt)\}_{\roots_0^+}$ by
the isomorpshims in Proposition \ref{prop:comparison of different dQEAs}.
We also define $\rsfs{U}_{q,X_{L_S\backslash G}}^k(\gt)$
as the inverse image of $\rsfs{U}_{q,X_{\roots}}^k(\gt)$ by
the immersion $X_{L_S\backslash G}(k)\longrightarrow X_{\roots}(k)$
defined in Proposition \ref{prop:embedding of the moduli}.
\end{defn}
\begin{rema}
By construction of the isomorphisms in Proposition \ref{prop:comparison of different dQEAs}, we can obtain an injection from $f^*\tilde{U_q^k(\lt)}$ to $\rsfs{U}_{q,X_{L_S\backslash G}}^k(\gt)$ where $\map{f}{X_{L_S\backslash G}(k)}{\spec k}$ is the canonical morphism.
\end{rema}
\begin{rema}
Since $X_{L_S\backslash G}(k)$ is an affine $k$-scheme 
and $\rsfs{U}_{q,X_{L_S\backslash G}}^k(\gt)$ is quasi-coherent,
we can reconstruct the sheaf from its global sections $U_{q,X_{L_S\backslash G}}^k(\gt) := \rsfs{U}_{q,X_{L_S\backslash G}}^k(\gt)(X_{L_S\backslash G}(k))$.

By construction of $\rsfs{U}_{q,X_{\roots}}^k(\gt)$, we have the local freeness of $U_{q,X_{L_S\backslash G}}^k(\gt)$. For $S \subset \simples$ with
$\abs{\simples \setminus S} = 1$, we can say a stronger statement:
there is a ``PBW basis'' for $U_{q,X_{L_S\backslash G}}^k(\gt)$,
constructed as follows.

At first one should note that $X_{L_S\backslash G}$ is covered by $U(\roots^+)$ and $U(\roots_0^+)$ where $\roots_0^+ = \roots_S^+\cup(\roots^-\setminus\roots_S^-)$. Take a reduced expression $s_{\bs{i}}$ of $w_0$ compatible with $\roots_0^+$. By definition of deformed quantum enveloping algebra, we have $\hE_{\bs{i}}$ and $\hK_{\lambda}$ as elements of $U_{q,X_{L_S\backslash G}}^k(\gt)$. For $\hF_{\bs{i},\alpha}$, we define
a global section by gluing $\psi_{\alpha}\hF_{\bs{i},\alpha}$ on $U(\roots^+)$ and $-\psi_{\alpha}\hF_{\bs{i},\alpha}$ on $U(\roots_0^+)$.
Then these global sections behave as ``quantum root vectors'' in $U_{q,X_{L_S\backslash G}}^k(\gt)$.

 For general $\frl_S \subset \g$, it seems to be still possible to show the freeness. For example we can define $F_{\alpha + \beta} \in U_{q,X_{H\backslash SL_3}}^k(\frak{sl}_3^{\sim})$ as follows:
\begin{align*}
 F_{\alpha + \beta} = \psi_{\alpha + \beta}F_{\beta}F_{\alpha}
-\psi_{\alpha + \beta}(q\psi_{\alpha} + q^{-1}\psi_{-\alpha})F_{\alpha}F_{\beta}.
\end{align*}
But we do not have any idea to obtain good
quantum root vectors in general case. Even in the case of $\frak{sl}_3$,
we do not know which element should be taken as $F_{\alpha + \beta}$.
\end{rema}

\subsection{Explicit construction of all quantizations}

Let $k$ be a $\Q$-algebra. In this case we can consider
$\rsfs{U}_{h,X_{\roots}}^k(\gt)$ and $\rsfs{U}_{h,X_{L_S\backslash G}}^k(\gt)$ on $X_{\roots}(k)$ and $X_{L_S\backslash G}(k)$
respectively. Also note that we can use $\phi_{\alpha} := 2\psi_{\alpha} - 1$ as generators of $\cl{O}(X_{L_S\backslash G}(k))$ since $2 \in k^{\times}$.

Let us recall the construction of the $2$-cocycle $\Psi_{h,\chi}^k$.
This can be characterized as a unique element inducing
the isomorphism (\ref{eq:associator for twisted induction}), which only depends on the $Q^-$-grading
of the $(\roots_0^+,\chi)$-twisted parabolically induced modules.
Since the identification in Proposition \ref{prop:comparison of different dQEAs}
preserves the triangular decomposition and the $Q$-grading,
we can see that $\Psi_{h,\chi}^k$ does not depend on the choice of $\roots_0^+$. This leads us to the following conclusion.
Note that there is a canonical homomorphism from $\cl{O}(X_{L_S\backslash G}(k))$ to $k$ when we take a character $\map{\chi}{2Q_0^-}{k}$ with $\roots_S^+ \subset \roots_0^+, \chi_{-2\alpha} = 1$ for $\alpha \in \roots_S^+$ and $\chi_{-2\alpha} -1 \in k^{\times}$ for $\alpha \in \roots_0^+\setminus\roots_S^+$.

\begin{thrm} \label{thrm:globalization}
There is a unique $2$-cocycle $\Psi_{h,\phi}^k \in U_h^{\cl{O}(X_{L_S \backslash G}(k))}(\g\times\g\times\frl_S)^{\frl_S}$
which coincides with $\Psi_{h,\chi}^k$ under the homomorphism above
for all $(k,\chi)$.
Moreover it has the following expansion:
\begin{align} \label{eq:expansion of the global 2-cocycle}
 \Psi_{h,\phi} = 1 - h\sum_{\alpha \in \roots^+\setminus\roots_S^+} d_{\alpha}(\phi_{\alpha} + 1)E_{\alpha}\tensor  F_{\alpha} \tensor 1 + \cdots.
\end{align}
\end{thrm}

As an immediate corollary, we can prove \cite[Proposition 5.3]{MR1817512} in a constructible way.

\begin{coro}
There is an algebraic family
$\{\cl{O}_{h,\phi}(L_S \backslash G)\}_{\phi \in X_{L_S \backslash G}}$ of equivariant deformation quantization of $\cl{O}(L_S \backslash G)$ with the associated
Poisson bracket $\{\tend,\tend\}_{\phi}$.
\end{coro}

By using \cite[Proposition 5.6]{MR1817512} we also obtain
the classification theorem. A \emph{formal path} on $X_{L_S \backslash G}$
is $(\phi_{\alpha}(h))_{\alpha \in \roots\setminus\roots_S}
\in \C\bbh^{\roots\setminus\roots_S}$ satisfying the relations
({\ref{relation:(i)}}), ({\ref{relation:(ii)}}), ({\ref{relation:(iii)}}). Then we can consider a homomorphism $\cl{O}(X_{L_S \backslash G})\longrightarrow \C\bbh$ by setting $\phi_{\alpha}\longmapsto \phi_{\alpha}(h)$
and obtain a deformation quantization $\cl{O}_{h,\phi(h)}(L_S \backslash G)$

\begin{coro} \label{coro:classification}
Let $\cl{O}_h(L_S \backslash G)$ be an equivariant deformation quantization
of $\cl{O}(L_S \backslash G)$. Then there is a unique formal path $\phi(h)$
on $X_{L_S \backslash G}$ such that $\cl{O}_h(L_S \backslash G)$ is equivalent
to $\cl{O}_{h,\phi(h)}(L_S \backslash G)$ as an equivariant deformation
quantization.
\end{coro}

\section{Comparison theorem}

Let $S \subset S_0 \subset \simples$ be subsets. The Levi subalgebra of $\g$ (resp. Levi subgroup of $G$)
associated to $S_0$ is denoted by $\frl_0$ (resp. $L_0$).
Then, by the results
in the previous section, we have the $2$-cocycle $\Psi_0 \in U_h^{\cl{O}(X_{L_S\backslash L_0})}(\frl_0\times \frl_0\times \frl_S)^{\frl_S}$.
Then we can thought $\Psi_0$ as a $2$-cocycle in $U_h^{\cl{O}(X_{L_S\backslash L_0})}(\g\times\g\times \frl_S)^{\frl_S}$ and obtain a deformation quantization
by the evaluation at an element of $X_{L_S\backslash L_0}$.
In light of Corollary \ref{coro:classification},
it is natural to compare this quantization and
quantizations coming from $\Psi \in U_h^{\cl{O}(X_{L_S \backslash G})}(\g\times\g\times\frl_S)^{\frl_S}$.
For such a purpose we consider a homomorphism
$\map{\pi}{\cl{O}(X_{L_S \backslash G}(k))}{\cl{O}(X_{L_S\backslash L_0}(k))}$ defined by
\begin{align*}
 \pi(\phi_{\alpha}) =
\begin{cases}
 \phi_{\alpha} & (\alpha \in \roots_{S_0}\setminus\roots_S), \\
 -1 & (\alpha \in \roots^+\setminus\roots_{S_0}^+), \\
1 & (\alpha \in \roots^-\setminus\roots_{S_0}^-).
\end{cases}
\end{align*}
\begin{prop} \label{prop:reduction of cocycle}
Let $\Psi \in U_h^{\cl{O}(X_{L_S \backslash G}(k))}(\g\times \g\times \frl_S)^{\frl}$
be the $2$-cocycle constructed in the previous section.
Then $\pi_*\Psi = \Psi_0$.
\end{prop}

We prove this fact as a corollary of the fully faithfulness of
a functor introduced in the following discussion, which holds in the integral setting.

Let $\roots_0^+$ be a positive system of $\roots$
containing $\roots_S^+\cup (\roots^-\setminus\roots_{S_0}^-)$.
Take a reduced expression $s_{\bs{i}}$ compatible with $\roots_0^+$.
Then we can introduce the notion of generalized $S_0$-maximal vector
in the similar manner to Definition \ref{defn:generalized maximal vector}.

We say that a $U_q^k(\g)$-module $V$ is
\emph{locally strongly $U_q^k(\n^+)$-finite} if,
for any $v \in V$, $U_q^k(\n^+)_{\mu}v = 0$ for all $\mu \in Q^+$ with finitely many exceptions.

\begin{thrm} \label{thrm:commutativity}
Let $(k,\chi)$ be an $\A[2Q_0^-]$-algebra with
$\chi_{-2\alpha} = 0$ for $\alpha \in \roots^-\setminus\roots_{S_0}^-$and $\chi_{-2\alpha} = 1$ for $\alpha \in \roots_S^+$.
\begin{enumerate}
 \item For any $V \in \mod{U_{q,\chi}^k}(\lot)$ and $W \in \mod{U_q^k(\g)}$, each $m \in W\tensor V$ has a unique generalized $S_0$-maximal vector $(m_{\Lambda})_{\Lambda}$ with $m_0 = m$.
 \item The functor $\map{\pdqind}{\mod{U_{q,\chi}^k(\lot)}}{\mod{U_{q,\chi}^k(\gt)}}, V\longmapsto U_{q,\chi}^k(\gt)\tensor_{U_{q,\chi}^k(\pot)} V$ is fully faithful.
 \item If a $U_{q,\chi}^k(\gt)$-module $W$ is locally strongly $U_q^k(\n^+)$-finite, there is a unique homomorphism $\pdqind{(W\tensor V)} \longrightarrow W\tensor\pdqind{V}$ which sends $m$ to a $S_0$-maximal vector
of the form $m_{(W)}\tensor (1\tensor m_{(V)}) + \cdots$. Moreover this is an isomorphism.
 \item Let $W, W'$ be locally strongly $U_q^k(\frb^+)$-finite $U_q^k(\g)$-module
and $V$ be a $U_{q,\chi}^k(\lot)$-module. Then the following
diagram commutes.
\begin{align*}
 \xymatrix{
\pdqind{(W\tensor W'\tensor V)}
\ar[r]^{\cong} \ar[rd]_{\cong} & W\tensor \pdqind{(W'\tensor V)} \ar[d]^{\cong}\\
& W\tensor W'\tensor \pdqind{V}
}
\end{align*}
\end{enumerate}

\end{thrm}


\begin{proof}
(i) The proof is very similar to that of Proposition \ref{prop:uniqueness of generalized maximal vector}, hence that of Lemma \ref{lemm:uniqueness and existence of partly maximal vector}.
The only different point is to show the upper diagonality and
the invertibility of diagonal entries in the usual sense.

Fix a multi-index $\Gamma$ and
consider its decomposition $\Gamma = \Gamma' + \delta_i$
such that $\gamma_1 = \gamma_2 = \cdot = \gamma_{i - 1} = 0$.
Then we look at the coefficient of $\hF_{\bs{i}}^{(\Gamma')}$
in (\ref{eq:recursion formula}) with $\beta = \alpha^{\bs{i}}_i$:
\begin{align*}
 \sum_{\Lambda} \hE_{\bs{i},\beta,(1)}m_{\Lambda,(W)}\tensor (\hE_{\bs{i},\beta,(2)}\hF_{\bs{i}}^{(\Lambda)})_{(\fru_S^-)}\tensor \pi_S((\hE_{\bs{i},\beta,(2)}\hF_{\bs{i}}^{(\Lambda)})_{(\frp_S)})m_{\Lambda,(V)}.
\end{align*}
By (\ref{eq:coproduct of E}), we can replace $\Delta(\hE_{\bs{i},\beta})$ by $K_{\Lambda^r\cdot\alpha^{\bs{i}}}^{-1}\tensor\hE_{\bs{i}}^{\Lambda^r}$ with $\Lambda^r \le i \le \Lambda^l$. Moreover
it suffices to consider the terms with $\Pi(\Lambda^l\cdot\alpha^{\bs{i}}) = 0$ and $\Lambda^r \neq \delta_i$ to show the upper diagonality.
Note that $\Lambda^r < i$ holds in this case.

Consider $\Lambda$ such that $\Lambda \in B_{l,\Pi(\Gamma\cdot\alpha^{\bs{i}})}$ and $\lambda_1 = \lambda_2 = \cdots = \lambda_{i - 1} = 0$.
Then we have the following expansion by Proposition \ref{prop:mixed LS relation} (i):
\begin{align*}
 \tE_{\bs{i}}^{\Lambda^r}\aF_{\bs{i}}^{(\Lambda)}
= \sum_{\substack{\Lambda^+ < i\le \Lambda^- \\ \Lambda^r\cdot\alpha^{\bs{i}} - \Lambda^+\cdot\alpha^{\bs{i}} = \Lambda\cdot\alpha^{\bs{i}} - \Lambda^-\cdot\alpha^{\bs{i}} \in Q^+}} C_{\Lambda^{\pm}} \aF_{\bs{i}}^{(\Lambda^-)}\tE_{\bs{i}}^{\Lambda^+}.
\end{align*}
Hence we also have the following:
\begin{align*}
 \hE_{\bs{i}}^{\Lambda^r}\hF_{\bs{i}}^{(\Lambda)}
= \sum_{\substack{\Lambda^+ < i\le \Lambda^- \\ \Lambda^r\cdot\alpha^{\bs{i}} - \Lambda^+\cdot\alpha^{\bs{i}} = \Lambda\cdot\alpha^{\bs{i}} - \Lambda^-\cdot\alpha^{\bs{i}} \in Q^+}} c_{\Lambda^{\pm}} \hF_{\bs{i}}^{(\Lambda^-)}\hE_{\bs{i}}^{\Lambda^+}. 
\end{align*}
If $\Lambda^+ \neq 0$, such a term is removed by $\pi_{S_0}$
since $\Lambda^+ < i$. It means that we only have to look at
the term with $\Lambda^+ = 0$ and $\Lambda^- = \Gamma'$ to
determine the coefficient of $\hF_{\bs{i}}^{(\Gamma')}$.
On the other hand, the coefficient of such a term can be calculated
as follows:
\begin{align*}
c_{0,\Gamma'} = C_{0,\Gamma'}\chi_{2\Lambda\cdot\alpha^{\bs{i}} - 2\Lambda^-\cdot\alpha^{\bs{i}}} = C_{0,\Gamma'}\chi_{2\Lambda^r\cdot\alpha^{\bs{i}}} = 0.
\end{align*}
Hence our $\Lambda$ does not affect the coefficient of $\hF_{\bs{i}}^{(\Gamma')}$. This shows the upper diagonality of $A_{l,\nu}$.

To show the invertibility, consider the case of $\Lambda = \Gamma$.
By the consideration above we may assume $\Lambda^r = \delta_i$.
Then a completely same argument in the proof of Lemma \ref{lemm:uniqueness and existence of partly maximal vector} works to prove
the invertibility of the diagonal entry. Then the upper triangularity
implies the invertibility of $A_{l,\nu}$.

(ii) This is a corollary of (i).

(iii) To see the former half of the statement,
it suffices to show that $m_{\Lambda} = 0$ except finitely many $\Lambda$ for any generalized $S_0$-maximal vector $(m_{\Lambda})_{\Lambda}$
since it implies the existence and uniqueness of required
$S_0$-maximal vectors and also the existence and uniqueness of the
required homomorphism.

By setting $W = U_q^k(\g)$ and $V = U_{q,\chi}^k(\lot)$, we can 
obtain the $S_0$-maximizer $(U_{\Lambda})_{\Lambda}$.
Moreover the consideration above implies that
$U_{\Lambda} \in \sum_{\Pi(\nu) = \Pi(\Lambda\cdot\alpha^{\bs{i}})}U_q^k(\frb^+)_{\nu}\tensor U_{q,\chi}^k(\lot)$. Hence our assumption on $W$ implies that
$m_{\Lambda} = U_{\Lambda}m_0 = 0$ except for finitely many $\Lambda$.

Next we show the bijectivity. Note that we have $\Delta(\hF_{\bs{i}}^{\Lambda}) = K_{\Lambda\cdot\alpha^{\bs{i}}}^{-1}\tensor \hF_{\bs{i}}^{\Lambda}$. Hence the image of $w\tensor (\hF_{\bs{i}}^{(\Gamma)}\tensor v)$ is 
\begin{align*}
 \sum_{\Lambda} K_{\Gamma\cdot\alpha^{\bs{i}}}^{-1}U_{\Lambda,(U_q^k(\frb^+))}w\tensor \hF_{\bs{i}}^{\Gamma}\hF_{\bs{i}}^{\Lambda}\tensor U_{\Lambda,(U_{q,\chi}^k(\frl_0{}^{\sim}))}v.
\end{align*}
Now we regard the homomorphism as an endomorphism $T$ on $U_{q,\chi}^k(\fru_{S_0}^-{}^{\sim}) \tensor W\tensor V$.
Then $U_{q,\chi}^k(\fru_{S_0}^-{}^{\sim})\tensor U_q^k(\frb^+)w \tensor V$ is stable under this map. Moreover our assumption implies the
existence of a positive integer $n$ such that $(T - \id)^n = 0$
on this module. This fact shows the bijectivity.

(iv) By (i) we have an automorphism $T$ on $W\tensor W'\tensor V$
such that $\pdqind{T}$ coincides with the following composition:
\begin{align*}
 \pdqind{(W\tensor W'\tensor V)}
&\longrightarrow W\tensor \pdqind{(W'\tensor V)} \\
&\longrightarrow W\tensor W'\tensor \pdqind{V}
\longrightarrow \pdqind{(W\tensor W'\tensor V)}.
\end{align*}
Then the image of
$1 \tensor (w\tensor w' \tensor v)$ at $W\tensor W'\tensor \pdqind{V}$
is of the form $T(w \tensor w'\tensor v)_{(W\tensor W')}\tensor (1\tensor T(w\tensor w'\tensor v)_{(V)}) + \cdots$. On the other hand,
by using the $S_0$-maximizer, we can calculate the image and
see that it is of the form $w\tensor w' \tensor (1\tensor v) + \cdots$. Hence we have $T = \id$ and the diagram commutes.
\end{proof}

\begin{proof}[Proof of Proposition \ref{prop:reduction of cocycle}]
By the construction of the cocycles, it suffices to show the commutativity of the following diagram:
\begin{align*}
 \xymatrix{
\dqind{(W\tensor V)} \ar[d] \ar[rr] && W\tensor \dqind{V}\\
\pdqind{\mathrm{ind}_{\frl_S,q}^{\frl_0,\chi}(W\tensor V)} \ar[r] & \pdqind{(W\tensor \mathrm{ind}_{\frl_S,q}^{\frl_0,\chi} V)} \ar[r] & W\tensor \pdqind{\mathrm{ind}_{\frl_S,q}^{\frl_0,\chi} V} \ar[u].
}
\end{align*}
Since the both images of $1 \tensor (w \tensor v)$ at the
upper right corner is of the form $w \tensor (1\tensor v) + \cdots$,
it suffices to show that these are $S_0$-maximal vectors.
The image under the upper isomorphism is trivially $S$-maximal.
To see the maximality of another image, we show more general statement: Let $W$ be a locally strongly $U_q^k(\n^+)$-finite module and
$V$ be a $U_{q,\chi}^k(\lot)$-module. If $m \in W\tensor V$
is $S$-maximal, the $S_0$-maximal vector $\tm = m_{(W)}\tensor (1\tensor m_{(V)}) + \cdots \in W\tensor \pdqind{V}$ is also $S$-maximal.

Take $\hE_{\bs{i},\beta}$ with $\beta \in \roots_{S_0}^+\setminus\roots_S^+$.
Then $\hE_{\bs{i},\beta}\tm$ is still of the form $(\hE_{\bs{i},\beta}m)_{(W)}\tensor (1\tensor (\hE_{\bs{i},\beta}m)_{(V)}) + \cdots$.
By a discussion similar to Lemma \ref{lemm:universal property}, $\hE_{\bs{i},\beta}\tm$ is again $S_0$-maximal. This implies
$\hE_{\bs{i},\beta}\tm = 0$ as a consequence of the uniqueness
for a maximal vector and $\hE_{\bs{i},\beta}m = 0$.
\end{proof}

At last we interpret this comparison theorem as a
statement on deformation quantizations. Let us recall the $U_h^k(\g)$-bimodule structure on $\cl{O}_h^k(G)$, induced from
the embedding $\cl{O}_h^k(G) \subset U_h^k(\g)^*$. Then,
for a $U_h^k(\frl_0)$-module algebra $B$,
we define the induced $U_h^k(\g)$-module algebra $\mathrm{ind}_{L_{0,h}}^{G_h} B$ as follows:
\begin{align*}
\mathrm{ind}_{L_{0,h}}^{G_h} B
= \{a \in \cl{O}_h^k(G)\tensor B\mid (r(a)\tensor \id)(x) = (\id\tensor l(x))(a)\text{ for all }x \in U_h^k(\frl_0)\},
\end{align*}
where $r(x)$ (resp. $l(x)$) is the right (resp. left) multiplication by $x \in U_h^k(\frl_0)$.
\begin{coro}
Let $(k\bbh,\chi)$ be an $\A[2Q_0^-]$-algebra same with $\chi_{-2\alpha} = 0$ for $\alpha \in \roots^-\setminus\roots_{S_0}^-$ and
$\chi_{-2\alpha} = 1$ for $\alpha \in \roots_S^+$.
Then $\cl{O}_{h,\chi}^k(L_S\backslash G)\cong \mathrm{ind}_{L_{0,h}}^{G_h} \cl{O}_{h,\chi}(L_S\backslash L_0)$.
\end{coro}
\begin{proof}
The $2$-cocycle arising from $U_h^k(\lt) \subset U_{h,\chi}^k(\frl_0{}^{\sim})$ (resp. $U_{h,\chi}(\gt)$) is denoted by $\Psi_0$ (resp. $\Psi$).
By Proposition \ref{prop:reduction of cocycle}, we have $\Psi_0 = \Psi$ throught the embedding $U_h^k(\frl_0\times\frl_0\times \frl_S) \subset U_h^k(\g\times\g\times\frl_S)$.

In light of the spectral decomposition (\ref{eq:spectral decomposition}) and the Peter-Weyl theorem for $\cl{O}_h^k(G)$, $\mathrm{ind}_{L_{0,h}}^{G_h} \cl{O}_{h,\chi}(L_S\backslash L_0)$ has the following expression:
\begin{align*}
&\mathrm{ind}_{L_{0,h}}^{G_h} \cl{O}_{h,\chi}(L_S\backslash L_0)
\cong \bigoplus_{\pi \in \irr_h\g, \rho \in \irr_h\frl_0} \Hom_{U_h^k(\frl_S)}(V_{\rho},k\bbh)\tensor \Hom_{U_h^k(\frl_0)}(V_{\pi},V_{\rho})\tensor V_{\pi}, \\
&(f\tensor T \tensor v)\ast (g\tensor S\tensor w)
= (f\tensor g)\Psi_0^{-1}\tensor (T\tensor S)\tensor (v\tensor w).
\end{align*}
Note that we have
\begin{align*}
(f\tensor g)\Psi_0^{-1}(T\tensor S) = (f\tensor g)(T\tensor S)\Psi^{-1} = (f\circ T\tensor g\circ S)\Psi^{-1}.
\end{align*}
On the other hand, we have
\begin{align*}
\Hom_{U_h^k(\frl_S)}(V_{\pi},k\bbh)
\cong
\bigoplus_{\rho \in \irr_h\frl_0} \Hom_{U_h^k(\frl_S)}(V_{\rho},k\bbh)\tensor \Hom_{U_h^k(\frl_0)}(V_{\pi},V_{\rho}),
\end{align*}
in which $f\tensor T$ corresponds to $f\circ T$.
Combining these facts we obtain the desired isomorphism.
\end{proof}
\vspace{10pt}
\noindent
{\bf Acknowlegements.}
The author is grateful to Yasuyuki Kawahigashi for comments on this paper and his invaluable supports. He also thanks Yuki Arano for encouraging him to write this paper.



\end{document}